\documentclass[12pt]{article}
\usepackage{amssymb}
\usepackage{amsmath}
\usepackage{amsthm}
\usepackage{color}
\usepackage{comment}
\usepackage{graphicx}
\begin{document}

\def\bu{{\bf u}}
\def\bv{{\bf v}}
\def\bz{{\bf z}}
\def\ee{\varepsilon}
\def\myproof{\noindent {\bf Proof.}\\}
\def\pdf{PDF}

\newcommand{\removableFootnote}[1]{}

\newtheorem{theorem}{Theorem}
\newtheorem{conjecture}[theorem]{Conjecture}
\newtheorem{lemma}[theorem]{Lemma}


\title{Stochastic Regular Grazing Bifurcations}
\author{
D.J.W.~Simpson$^{\dagger}$, S.J.~Hogan$^{\ddagger}$ and R.~Kuske$^{\dagger}$\\\\
$^{\dagger}$Department of Mathematics\\
University of British Columbia\\
Vancouver, BC, V6T1Z2, Canada\\\\
$^{\ddagger}$Department of Engineering Mathematics\\
University of Bristol\\
Bristol, BS8 1TR, UK
}
\maketitle

\begin{abstract}

A grazing bifurcation corresponds to the
collision of a periodic orbit with a switching manifold
in a piecewise-smooth ODE system
and often generates complicated dynamics.
The lowest order terms of the induced Poincar\'{e} map
expanded about a regular grazing bifurcation constitute a Nordmark map.
In this paper we study a normal form of the Nordmark map in two dimensions
with additive Gaussian noise of amplitude, $\varepsilon$.
We show that this particular noise formulation arises in a general setting and
consider a harmonically forced linear oscillator
subject to compliant impacts to illustrate the accuracy of the map.
Numerically computed invariant densities of the stochastic Nordmark map
can take highly irregular forms, or,
if there exists an attracting period-$n$ solution when $\varepsilon = 0$,
be well approximated by the sum of $n$ Gaussian densities
centred about each point of the deterministic solution, and scaled by $\frac{1}{n}$,
for sufficiently small $\varepsilon > 0$.
We explain the irregular forms and calculate
the covariance matrices associated with the Gaussian approximations
in terms of the parameters of the map.
Close to the grazing bifurcation
the size of the invariant density may be proportional to $\sqrt{\varepsilon}$,
as a consequence of a square-root singularity in the map.
Sequences of transitions between different dynamical regimes
that occur as the primary bifurcation parameter is varied
have not been described previously.

\end{abstract}

\section{Introduction}
\label{sec:INTRO}
\setcounter{equation}{0}

\begin{figure}[b!]
\begin{center}
\setlength{\unitlength}{1cm}
\begin{picture}(14.67,3.8)
\put(0,0){\includegraphics[height=3.5cm]{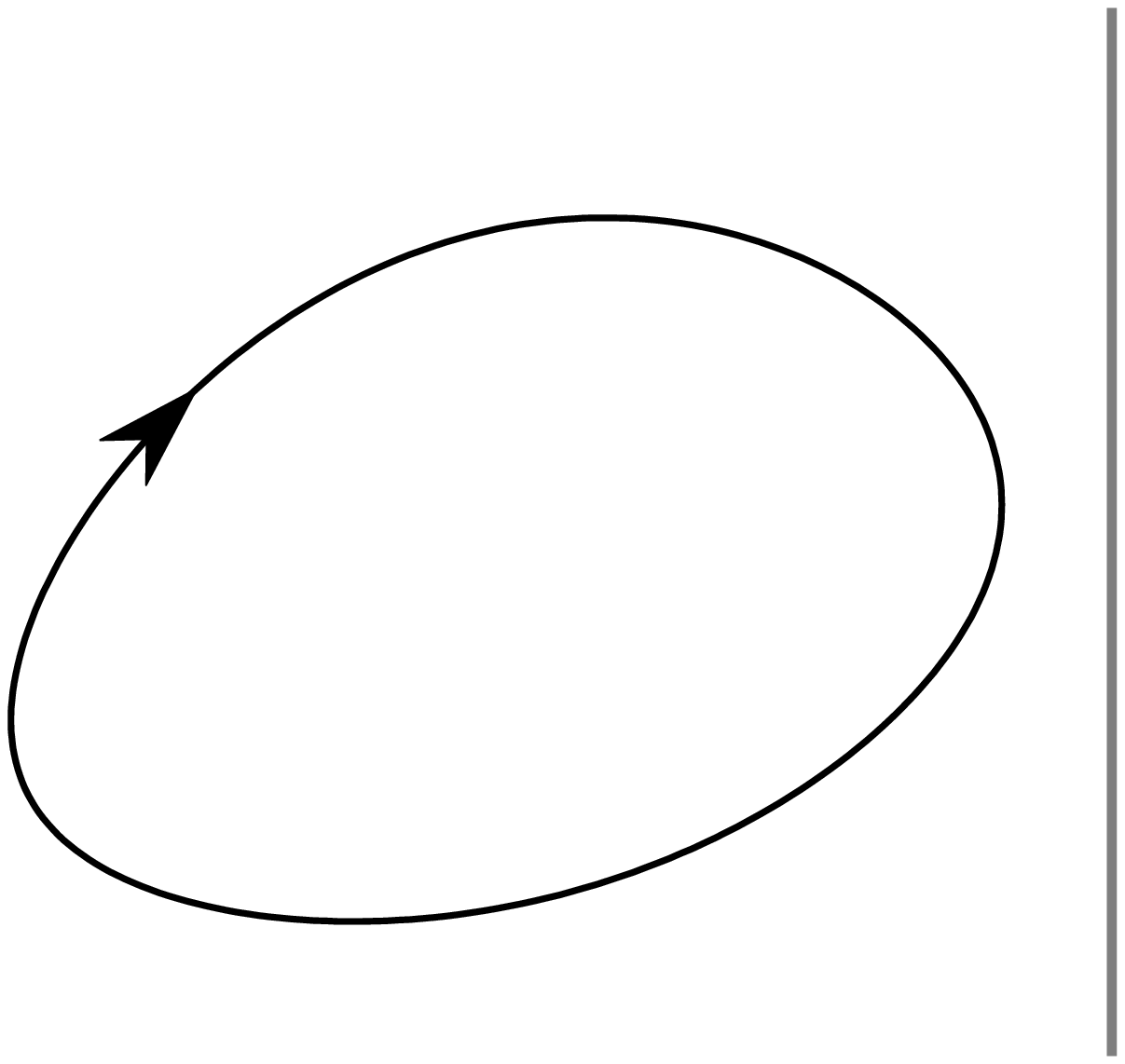}}
\put(5,0){\includegraphics[height=3.5cm]{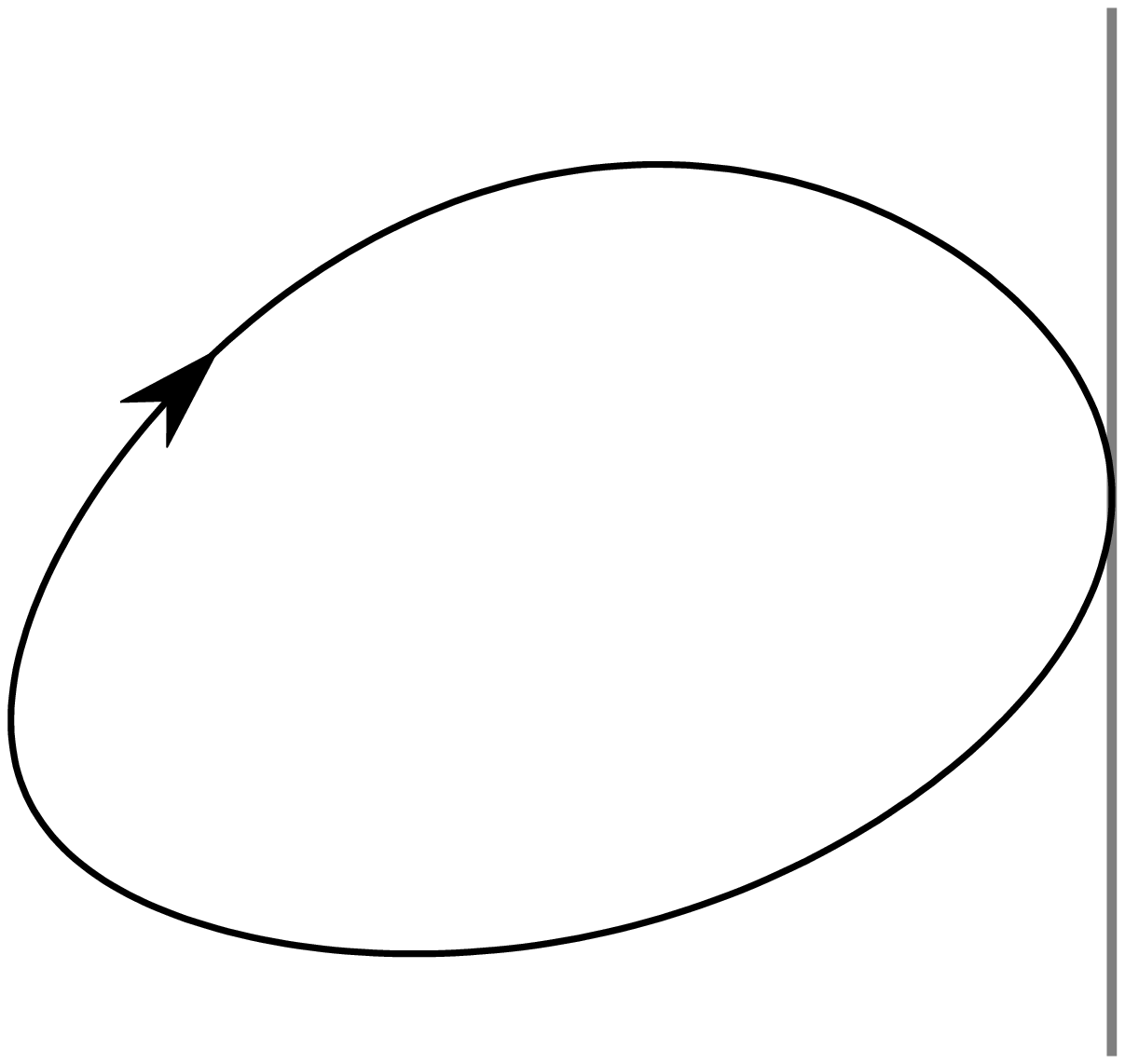}}
\put(10,0){\includegraphics[height=3.5cm]{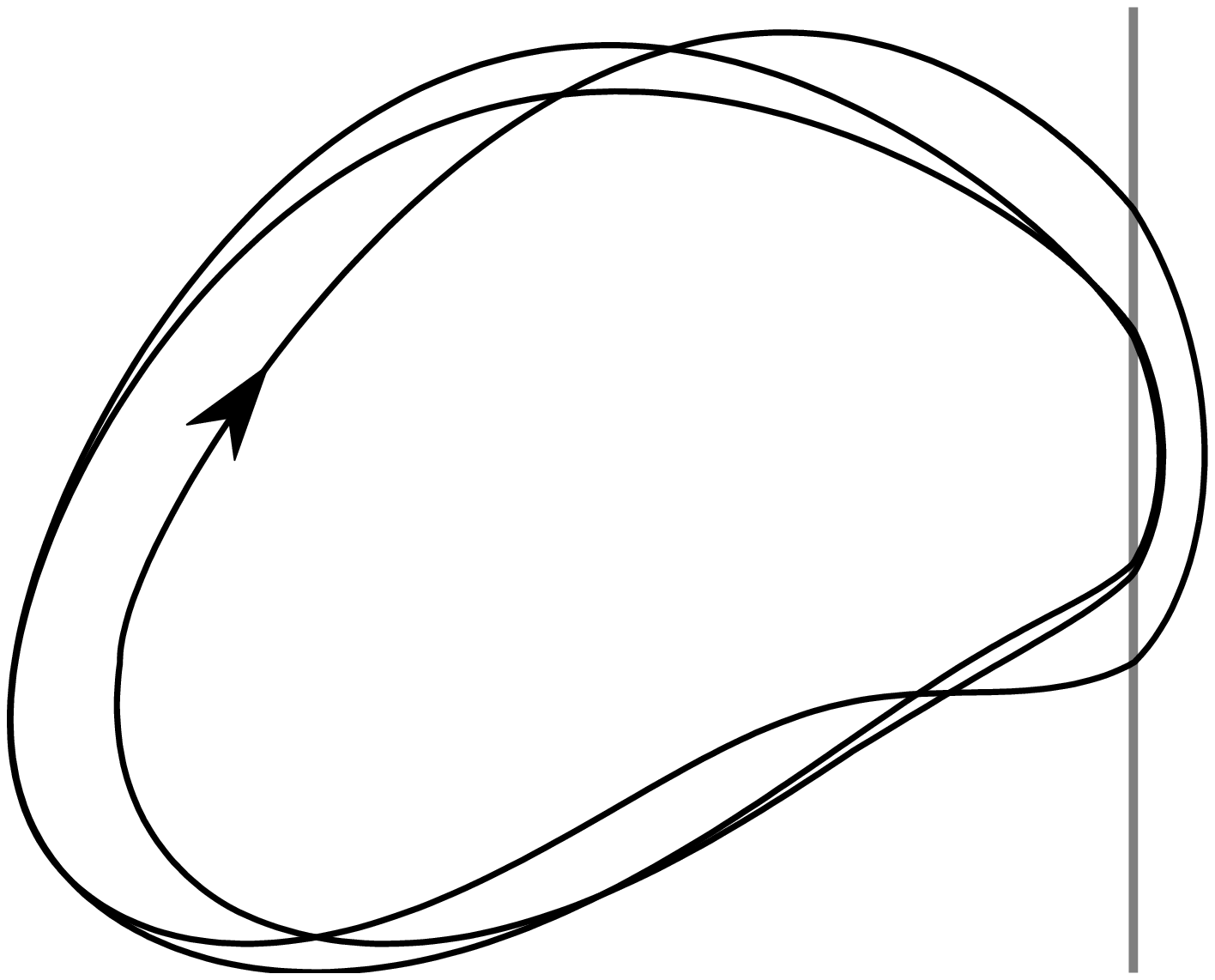}}
\put(2,3.7){\small $\mu<0$}
\put(7,3.7){\small $\mu=0$}
\put(12,3.7){\small $\mu>0$}
\put(4.4,.3){\small $\Sigma$}
\put(9.4,.3){\small $\Sigma$}
\put(14.4,.3){\small $\Sigma$}
\end{picture}
\caption{
Phase space schematics illustrating a grazing bifurcation.
The vertical line, $\Sigma$, represents a switching manifold.
The bifurcation occurs when a periodic orbit intersects the switching manifold
as a parameter, $\mu$, is varied.
Dynamics shortly after the bifurcation may be extremely complex
and depend on the stability multipliers of the
non-impacting periodic orbit and the precise nature of the grazing scenario \cite{DiBu08}.
\label{fig:grazBifSchem}
}
\end{center}
\end{figure}

Deterministic, piecewise-smooth systems are commonly used
as mathematical models of vibro-impacting systems \cite{WiDe00,Br99,BlCz99,Ib09}.
Specific examples include
atomic force microscopy \cite{DaZh07,RaMe08},
gear assemblies \cite{ThNa00,HaWi07},
metal cutters \cite{Gr88,Wi97},
and vibrating heat-exchanger tubes \cite{PaLi92,DeFr99}.
Piecewise-smooth systems are characterized by the presence of
codimension-one regions of phase space, termed switching manifolds,
across which the functional form of the system changes.
In the context of vibro-impacting systems, switching manifolds
correspond to an impact or loss of contact.
Motion that involves recurring impacts is often born
in the collision of an attracting non-impacting periodic orbit
with a switching manifold. 
This collision is known as a grazing bifurcation, Fig.~\ref{fig:grazBifSchem}.
Grazing bifurcations have been intensely analyzed 
and may be the cause of complex dynamics including chaos,
see \cite{DiBu08} and references within.
In this paper we study the effect of noise on
a common class of grazing bifurcations by analyzing a stochastic version
of the induced return map.

The effect of small noise on a well-understood deterministic system is often
intuitive and qualitatively predictable
but in certain situations may produce interesting new dynamics.
For instance there are various mechanisms,
such as coherence resonance, by which
the addition of noise to quiescent systems
induces relatively regular oscillations \cite{LiGa04,PiKu97,MuVa05},
as well as more complicated dynamics such as mixed-mode oscillations \cite{DeGu12}. 
Alternatively noise may suppress chaos \cite{Ra95}.
Studies of stochastic versions of some prototypical models
are found in \cite{HoLe06,Da98,ArBl99,ZhLu93,Sc96}.

In \cite{GrHo12,Gr05}, Griffin and Hogan
study a one-dimensional piecewise-linear map with additive noise.
For uniformly distributed noise
they calculate basins of attraction and widths of invariant densities.
With Gaussian noise the map resists such a precise analysis
and instead numerical results are presented revealing that the two types of noise
produce qualitatively similar behaviour.
In a period-incrementing scenario,
the invariant density of the stochastic map may undergo several transitions
(as distinguished by the number of peaks in the density)
as the noise amplitude is varied \cite{GrHo05}.
Numerical simulations of a two-dimensional ODE system with additive noise
designed such that the relevant Poincar\'{e} map is piecewise-linear,
give essentially the same dynamics as the map.
In \cite{GrHo12b},
numerical simulations of a noisy piecewise-smooth map are shown to
exhibit many similarities to the output of a DC/DC converter.
Noise-induced stabilization in one-dimensional discontinuous maps
is described in \cite{Wa98}.

Impact oscillators subject to noise have been studied
from the point of view of characterizing
small-amplitude, noise-induced oscillations.
The early work of Dimentberg and Menyailov \cite{DiMe79} concerns
an unforced, damped, linear oscillator experiencing noise and
instantaneous impacts with either one wall or two symmetrically placed walls.
Zhuravlev's transformation \cite{Zh76,Ib09}
converts the stochastic equation of motion with an impact rule
into a single discontinuous stochastic differential equation\removableFootnote{
If impacts are elastic, 
the invariant probability density function
for the position and velocity of the oscillator
is a truncated Gaussian.
Moreover the shape of the Gaussian is the same as that for the corresponding non-impacting problem
and indeed the theorem given in \cite{DiIo04} indicates this.
}.
To then determine the invariant density of the position and velocity of the oscillator,
in the case that impacts are slightly inelastic
the system may be treated as a perturbation from a Hamiltonian system.
Here stochastic averaging \cite{RoSp86,Kh66,AnAs02}
leads to a simple stochastic differential equation for the energy.
The work \cite{DiMe79} has been extended to
compliant impacts \cite{FoBr96},
Hertzian contacts \cite{DiIo04},
and nonlinear oscillators \cite{SrPa05}\removableFootnote{
The PDF for the position of the oscillator has been studied numerically \cite{IoSo06}.
For two symmetrically placed walls,
this PDF has peaks near the walls when the coefficient of restitution is small,
$r < 0.5$, say,
indicating that the system prefers to be near the walls in this case.
Also in \cite{FeXu09b} the authors consider
a single wall placed at the equilibrium position of the oscillator.
Here the PDF has a local minimum at the origin indicating that
noise induces a degree of oscillatory motion.
}.
The presence of forcing adds significant difficulty to the problem
but has recently been investigated numerically
and via stochastic averaging \cite{DiIo04,DiIo05}
and analyzed with series expansions and mean Poincar\'{e} maps
in the context of gear rattling \cite{FePf98,PfKu90}.

Grazing bifurcations are typically studied
by deriving and analyzing a Poincar\'{e} map
that is valid in a neighbourhood of the bifurcation.
Different grazing scenarios yield fundamentally different Poincar\'{e} maps.
In two dimensions, a normal form of the {\em Nordmark map} is:
\begin{equation}
\left[ \begin{array}{c} x' \\ y' \end{array} \right] =
\left\{ \begin{array}{lc}
A \left[ \begin{array}{c} x \\ y \end{array} \right] +
\left[ \begin{array}{c} 0 \\ 1 \end{array} \right] \mu \;, & x \le 0 \\
A \left[ \begin{array}{c} x \\ y - \chi \sqrt{x} \end{array} \right] +
\left[ \begin{array}{c} 0 \\ 1 \end{array} \right] \mu \;, & x \ge 0
\end{array} \right. \;,
\label{eq:NordmarkMap}
\end{equation}
where
\begin{equation}
A = \left[ \begin{array}{cc}
\tau & 1 \\ -\delta & 0
\end{array} \right] \;, \qquad
\chi = \pm 1 \;,
\end{equation}
and $\tau,\delta \in \mathbb{R}$.
This map applies to both {\em regular grazing} in
three-dimensional, piecewise-smooth systems \cite{DiBu01},
and grazing in vibro-impacting systems for which impacts
are modelled as instantaneous events with velocity reversal and energy loss \cite{No91,No97}.
The grazing bifurcation occurs at the origin when $\mu = 0$.
Since (\ref{eq:NordmarkMap}) is a truncated series expansion centred about the bifurcation,
(\ref{eq:NordmarkMap}) is primarily of interest for small values of $x$, $y$ and $\mu$.
The construction of (\ref{eq:NordmarkMap}) is described in \S\ref{sub:CONSTR}.
Here we note that the fixed point of the left half-map of (\ref{eq:NordmarkMap}) is given by
\begin{equation}
\left[ \begin{array}{c}
x^{*(L)} \\ y^{*(L)}
\end{array} \right]
= \frac{1}{1-\tau+\delta}
\left[ \begin{array}{c}
1 \\ 1-\tau
\end{array} \right] \mu \;,
\label{eq:xyStarL}
\end{equation}
and corresponds to the non-impacting periodic orbit.
Assuming this periodic orbit is attracting,
the fixed point, (\ref{eq:xyStarL}),
is admissible (i.e.~$x^{*(L)} \le 0$) when $\mu \le 0$.

Here we describe the basic results and summarize the remainder of the paper.
We analyze (\ref{eq:NordmarkMap}) with
small-amplitude, additive Gaussian noise:
\begin{equation}
\left[ \begin{array}{c} x' \\ y' \end{array} \right] =
\left\{ \begin{array}{lc}
A \left[ \begin{array}{c} x \\ y \end{array} \right] +
\left[ \begin{array}{c} 0 \\ 1 \end{array} \right] \mu \;, & x \le 0 \\
A \left[ \begin{array}{c} x \\ y - \chi \sqrt{x} \end{array} \right] +
\left[ \begin{array}{c} 0 \\ 1 \end{array} \right] \mu \;, & x \ge 0
\end{array} \right\} + \ee \xi \;,
\label{eq:P}
\end{equation}
where $0 < \ee \ll 1$, and
\begin{equation}
\xi \sim N(0,\Theta) \;,
\label{eq:xi}
\end{equation}
that is, $\xi$ is a two-dimensional Gaussian random variable with zero mean
and covariance matrix $\Theta$.
In \S\ref{sub:ADD} we demonstrate that this noise formulation
arises when small-amplitude Brownian motion is added to
the underlying differential equations
and obtain an explicit expression for $\Theta$.
For a model of a linear oscillator experiencing impacts with a compliant support,
in \S\ref{sub:IMPACT} we show that (\ref{eq:P})
exhibits quantitatively the same dynamics as the
stochastic differential equation model.
Calculations for this section are deferred to Appendix \ref{sec:APP}.

\begin{figure}[b!]
\begin{center}
\setlength{\unitlength}{1cm}
\begin{picture}(12.3,9.3)
\put(.3,.3){\includegraphics[height=9cm]{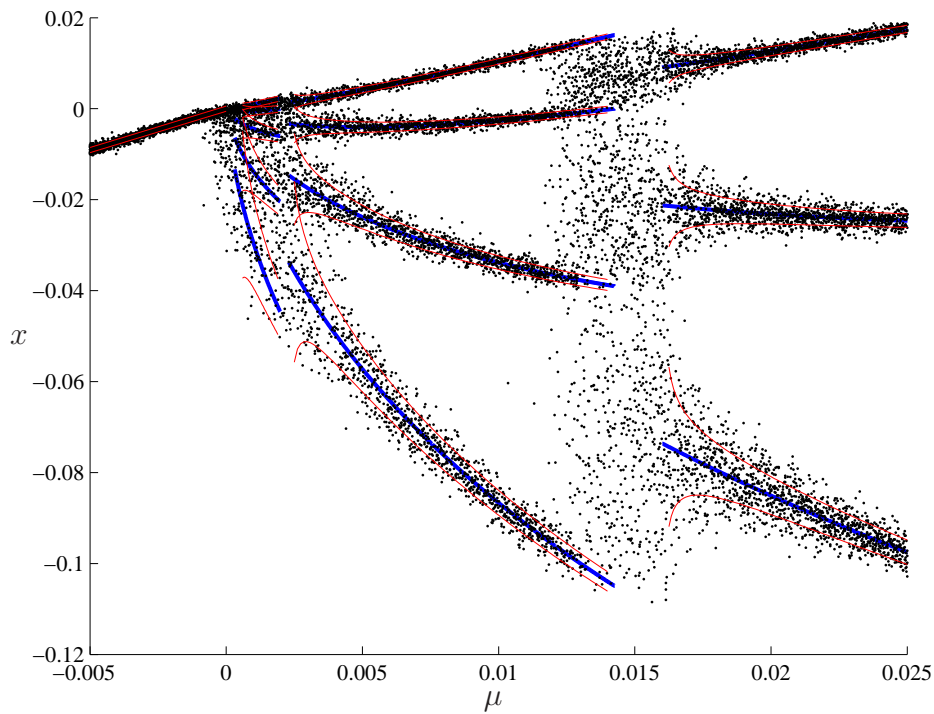}}
\put(6.3,0){$\mu$}
\put(0,4.8){$x$}
\end{picture}
\caption{
A bifurcation diagram of the stochastic Nordmark map, (\ref{eq:P}),
with $(\tau,\delta,\chi) = (0.5,0.05,1)$,
$\ee = 0.00025$,
and $\Theta = I$ (the $2 \times 2$ identity matrix).
The dots are numerically computed iterates of (\ref{eq:P})
with transients omitted.
The blue curves denote deterministic, admissible, attracting periodic solutions,
computed analytically in \S\ref{sub:MPS}.
The red curves approximate one standard deviation of the invariant density
from the blue curves, as obtained via linearization, see \S\ref{sub:GAUSS}.
\label{fig:stochBifDiag_c}
}
\end{center}
\end{figure}

\begin{figure}[b!]
\begin{center}
\setlength{\unitlength}{1cm}
\begin{picture}(12.3,9.3)
\put(.3,.3){\includegraphics[height=9cm]{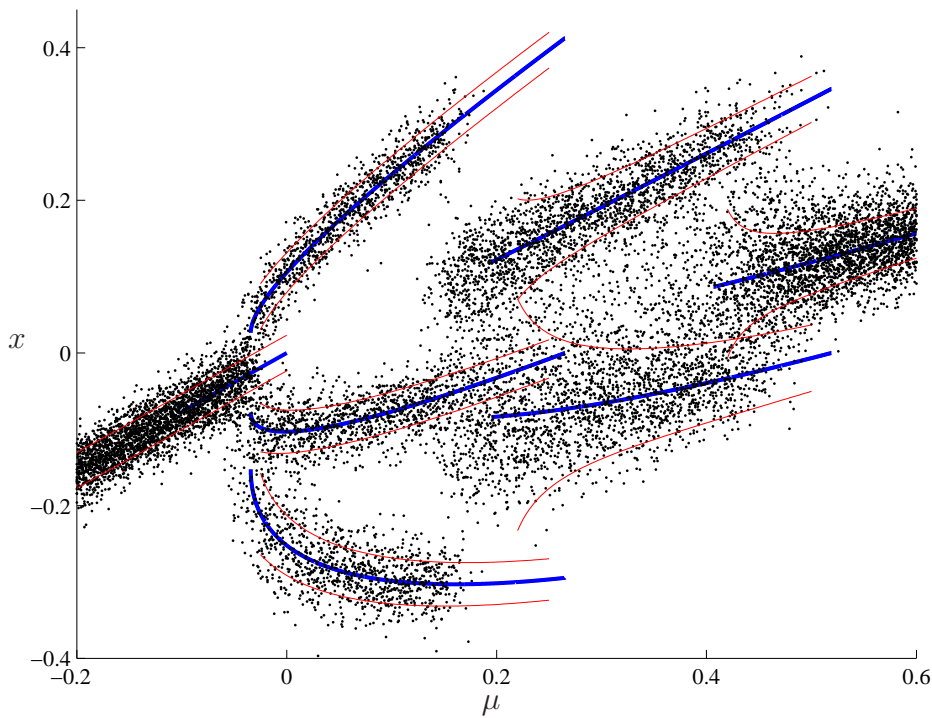}}
\put(6.3,0){$\mu$}
\put(0,4.8){$x$}
\end{picture}
\caption{
A bifurcation diagram of (\ref{eq:P}),
with $(\tau,\delta,\chi) = (0.2,0.5,1)$,
$\ee = 0.01$, and $\Theta = I$.
The meaning of various figure elements are explained
in the caption of Fig.~\ref{fig:stochBifDiag_c}.
\label{fig:stochBifDiag_b}
}
\end{center}
\end{figure}

Figs.~\ref{fig:stochBifDiag_c} and \ref{fig:stochBifDiag_b}
show bifurcation diagrams of (\ref{eq:P}) for two different
choices of the parameter values.
These were computed numerically
by iterating (\ref{eq:P}) over a range of fixed values of $\mu$.
As described in \cite{Gr05},
the noise blurs the underlying deterministic bifurcation diagram and
causes impacts to occasionally occur when $\mu < 0$\removableFootnote{
Also Griffin \cite{Gr05} observes that
a plot of the mean value of a numerical approximation to the invariant density
against $\mu$ is smoothed by the noise.
}.
In \S\ref{sub:MPS} we analytically compute periodic solutions of (\ref{eq:P}) with $\ee = 0$.
These are indicated by blue curves
in Figs.~\ref{fig:stochBifDiag_c} and \ref{fig:stochBifDiag_b}.
For the parameter values of Fig.~\ref{fig:stochBifDiag_c}, in the absence of noise
the maps exhibits a period-incrementing cascade, \cite{No01,DiBu08}.
For every $n \ge 1$, there exists an interval of positive $\mu$-values
for which period-$n$ solutions are admissible and attracting
(periods 3, 4 and 5 are indicated in Fig.~\ref{fig:stochBifDiag_c}).
These intervals are ordered by $n$,
disjoint, and separated by small chaotic bands.
As $n \to \infty$, the intervals approach $\mu = 0$ from the right.

When $\ee > 0$, dynamics may be roughly periodic
(e.g.~when $\mu = 0.01$ in Fig.~\ref{fig:stochBifDiag_c}),
however, the noise smothers high-period solutions. 
As with the logistic map \cite{CrFa82,MaHa81},
the highest period that is in some sense distinguishable,
decreases with an increase in the noise amplitude.
Figs.~\ref{fig:stochBifDiag_c} and \ref{fig:stochBifDiag_b}
suggest the existence of a stable invariant density
in a neighbourhood of the underlying attractor, for any $\mu$.
Throughout the paper we study this density
but do not formally prove its existence or uniqueness.
The invariant density may be approximately Gaussian about each point
of a deterministic periodic solution of relatively low period.
In \S\ref{sub:GAUSS} we derive the covariance matrices of these Gaussians.
(The Gaussian approximation to the unimodal invariant density
centred about the fixed point, (\ref{eq:xyStarL}), for $\mu < 0$,
may be computed in the same way.)
In \S\ref{sub:APPROX} we derive an approximation to these covariance matrices
and in \S\ref{sub:THETA} investigate the effect of different $\Theta$, (\ref{eq:xi}).
The near-Gaussian densities are discussed further in \S\ref{sub:D1}.
The accuracy of the Gaussian approximation decreases
as $\ee$ is increased, \S\ref{sub:D2}.
For relatively large $\ee$,
the density is strongly non-Gaussian and may
exhibit some interesting novel characteristics, \S\ref{sub:D3}.
Finally, we note that when $\ee = 0$ attracting solutions may coexist.
For instance in Fig.~\ref{fig:stochBifDiag_b},
for small $\mu < 0$ the fixed point (\ref{eq:xyStarL})
coexists with a period-3 solution,
and for a relatively wide range $\mu$-values in the middle of the figure,
this period-3 solution coexists with either a period-2 solution
or what we numerically observe to be a chaotic attractor\removableFootnote{
The period-2 solution loses stability in a subcritical period-doubling
bifurcation at $\mu \approx 0.19515$.
The saddle-type period-4 solution created in this bifurcation exists
to the right of this bifurcation.
Numerics suggest that at $\mu \approx 0.1963$ a 4-piece chaotic attractor
is created that persists to the left of this bifurcation.
The bifurcation here seems to be some sort of catastrophic discontinuity-induced bifurcation
at which some pieces of the basin of attraction of the 2-cycle break off
and form a basin for a 4-piece chaotic attractor.
The resulting boundary between the 2-cycle and the chaotic attractor
is formed by the stable and unstable manifolds of the 4-cycle.
As $\mu$ is decreased, the 4-piece chaotic attractor appears to merge into a 2-piece chaotic
attractor, then into a 1-piece chaotic attractor, before finally being destroyed
in some sort of catastrophe at $\mu \approx 0.093$.
}.
When $\ee > 0$, numerics indicate the existence of
a single invariant density with peaks at the deterministic attractors, \S\ref{sub:D4}.
Finally conclusions are presented in \S\ref{sec:CONC}.

\section{Derivation of the stochastic Nordmark map}
\label{sec:DERIV}
\setcounter{equation}{0}

The purpose of this section is to justify the manner by which
we have incorporated noise into the Nordmark map, (\ref{eq:P}) with (\ref{eq:xi}),
and derive an explicit expression for $\Theta$
in terms of an underlying stochastic differential equation.
Certainly noise may enter a vibro-impacting system in a variety of different ways.
Noise in external forcing may be represented by a deterministic signal
plus a weak random component modelled by white or coloured noise.
Model parameters may be treated stochastically.
Also impact events may be modelled randomly
(indeed experiments reveal impacts may induce high frequency oscillations that quickly decay
and vanish \cite{OeHi97}).
The manner by which randomness is modelled affects
the nature of the induced Poincar\'{e} map.
Below we consider additive, white Gaussian noise and derive (\ref{eq:P}).

\subsection{Construction of the deterministic Nordmark map}
\label{sub:CONSTR}

To describe the construction of the Nordmark map
it is necessary to introduce a general system that exhibits a regular grazing bifurcation.
To achieve this in a general setting
requires considerable space in order to carefully state assumptions,
non-degeneracy conditions and explain their meaning.
For brevity we omit most of these statements -- details
are given in \cite{DiBu08}.

Consider an ODE system in $\mathbb{R}^3$, $\bu = (u,v,w)$,
with a single smooth switching manifold,
\begin{equation}
\Sigma = \{ \bu ~|~ u=0 \} \;,
\label{eq:Sigma}
\end{equation}
that we write as
\begin{equation}
\dot{\bu} = \left\{ \begin{array}{lc}
f^{(L)}(\bu;\eta) \;, & u < 0 \\
f^{(R)}(\bu;\eta) \;, & u > 0
\end{array} \right. \;.
\label{eq:ODE}
\end{equation}
where $f^{(L)}$ and $f^{(R)}$ are smooth functions
extendable beyond their respective half-spaces,
and $\eta$ is a parameter.
Suppose when $\eta = 0$
there exists an attracting periodic orbit that intersects $\bu = 0$
but is otherwise contained in $\{ \bu ~|~ u < 0 \}$.
To negate the possibility of sliding motion \cite{DiBu08},
and specify a direction of flow, suppose
\begin{equation}
{\rm sgn} \left( e_1^{\sf T} f^{(L)}(0,v,w;\eta) \right) 
= {\rm sgn} \left( e_1^{\sf T} f^{(R)}(0,v,w;\eta) \right)
= {\rm sgn}(v) \;,
\label{eq:noSliding}
\end{equation}
in a neighbourhood of $\bu = 0$.

\begin{figure}[b!]
\begin{center}
\setlength{\unitlength}{1cm}
\begin{picture}(13.25,10)
\put(0,0){\includegraphics[height=10cm]{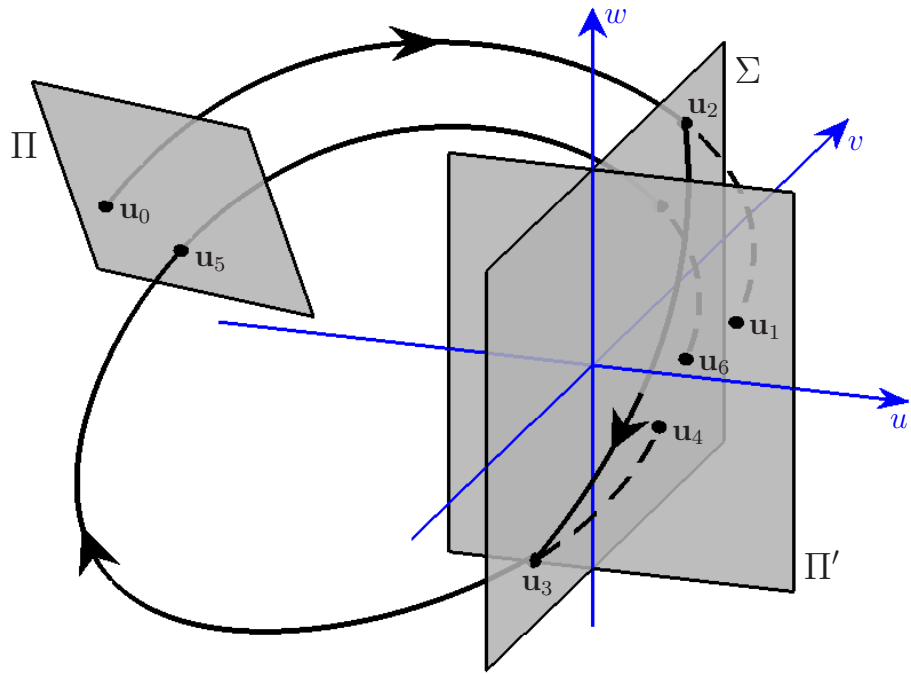}}
\put(12.8,4.1){\color{blue} $u$}
\put(12.24,7.8){\color{blue} $v$}
\put(9,9.5){\color{blue} $w$}
\put(1.1,7.7){\large $\Pi$}
\put(11.65,2.1){\large $\Pi'$}
\put(10.73,8.7){\large $\Sigma$}
\put(2.55,6.9){$\bu_0$}
\put(10.92,5.35){$\bu_1$}
\put(10.1,8.3){$\bu_2$}
\put(7.9,1.95){$\bu_3$}
\put(9.9,3.95){$\bu_4$}
\put(3.55,6.26){$\bu_5$}
\put(10.25,4.88){$\bu_6$}
\end{picture}
\caption{
A schematic for regular grazing in three dimensions.
We show intersections of a trajectory of the system (\ref{eq:ODE})
with a generic Poincar\'{e} section, $\Pi$,
the alternate Poincar\'{e} section, $\Pi'$, (\ref{eq:Pip}),
and the switching manifold, $\Sigma$, (\ref{eq:Sigma}).
Since the trajectory intersects the switching manifold at $\bu_2$,
the points $\bu_1, \bu_4 \in \Pi'$ are found by solving the left half-system, $f^{(L)}$,
in the right half-space.
The discontinuity map is $\bu_4 = D(\bu_1)$, (\ref{eq:D}).
The global map is $\bu_6 = G(\bu_4)$, (\ref{eq:G}).
The normal form Nordmark map, (\ref{eq:NordmarkMap}), is given by $G \circ D$,
after the coordinate change (\ref{eq:coordinateChange}) and truncation.
\label{fig:overallSchem}
}
\end{center}
\end{figure}

We introduce two Poincar\'{e} sections:
$\Pi$ is a generic cross-section transversal to the flow lying 
entirely in the left half-space, and 
\begin{equation}
\Pi' = \{ \bu ~|~ v=0 \} \;,
\label{eq:Pip}
\end{equation}
see Fig.~\ref{fig:overallSchem}.
To analyze oscillatory dynamics
one can study the induced return map on $\Pi$,
or a return map on $\Pi'$ using intersection points that for $u>0$ are virtual,
as explained below.
The Poincar\'{e} map on $\Pi$ is useful for comparing with numerical simulations.
However, we use the section $\Pi'$, because,
although some iterates are virtual, the map takes a simpler form
and is more amenable to analysis.

To obtain the return map on $\Pi'$,
consider the trajectory from a point, $\bu_0$, on $\Pi$, as governed by (\ref{eq:ODE}),
and suppose that this trajectory intersects $\Sigma$ at $\bu_2$ and $\bu_3$
before returning to $\Pi$ at some $\bu_5$, as shown in Fig.~\ref{fig:overallSchem}.
We use the left half-flow, $f^{(L)}$, to follow the trajectory
beyond $\bu_2$ into the right half-space until intersection with $\Pi'$ at $\bu_1$.
Similarly we use $f^{(L)}$ to follow the trajectory backwards from $\bu_3$ to $\bu_4$.
The map from $\bu_1$ to $\bu_4$ defined in this manner, call it $D$, is known as a
{\em discontinuity map} and takes the form
\begin{equation}
\left[ \begin{array}{c}
u_4 \\ w_4
\end{array} \right] =
D(u_1,w_1;\eta) =
\left\{ \begin{array}{lc}
\left[ \begin{array}{c}
u_1 \\ w_1
\end{array} \right] \;, & u_1 \le 0 \\
\left[ \begin{array}{c}
u_1 + O(u_1^{\frac{3}{2}}) \\
w_1 - (c + O(|w_1,\eta|^1)) \sqrt{u_1} + O(u_1)
\end{array} \right] \;, & u_1 \ge 0
\end{array} \right. \;,
\label{eq:D}
\end{equation}
for some scalar constant $c$\removableFootnote{
If we expand the half-flows about $(\bu,\eta) = (0,0)$ as:
\begin{equation}
\left[ \begin{array}{c}
\dot{u} \\ \dot{v} \\ \dot{w}
\end{array} \right] =
\left[ \begin{array}{c}
\alpha_J u + \beta_J v + O(2) \\ -\gamma_J + O(1) \\ \zeta_J + O(1)
\end{array} \right] \;,
\end{equation}
for $J = L,R$, then
\begin{equation}
c = \frac{2 \sqrt{2 \gamma_L}}{\sqrt{\beta_L}}
\left( \frac{\zeta_L}{\gamma_L} - \frac{\zeta_R}{\gamma_R} \right) \;.
\end{equation}
Note that $\beta_L, \beta_R > 0$ by (\ref{eq:noSliding})
(specifying the direction of flow across $\Sigma$).
Also, $\gamma_L, \gamma_R > 0$, as implied by Fig.~\ref{fig:overallSchem}
(specifying the direction of flow across $\Pi'$).
}.
This concept was originally introduced by Nordmark \cite{No91,No97}
in regards to instantaneous impacts and is nowadays a fundamental tool
for analyzing a variety of grazing phenomena
\cite{DiBu01,DaNo00,DiBu01c,FrNo97,FrNo00}.
The discontinuity map is nonsmooth and for $u_1>0$
corresponds to motion on an $O(\sqrt{u_1})$ time scale,
where $u_1$ is assumed to be small.
In contrast, the induced Poincar\'{e} map from $\bu_4$ to $\bu_6$,
as governed by the left half-flow, captures the global dynamics on an $O(1)$ time scale.
This global map, which we call $G$, is smooth and thus may be written as
\begin{equation}
\left[ \begin{array}{c}
u_6 \\ w_6
\end{array} \right] =
G(u_4,w_4;\eta) =
\hat{A} \left[ \begin{array}{c}
u_4 \\ w_4
\end{array} \right]
+ \hat{b} \eta + O(|u_4,w_4,\eta|^2) \;.
\label{eq:G}
\end{equation}
for some constant matrix $\hat{A}$ and vector $\hat{b}$.
The Nordmark map (\ref{eq:NordmarkMap}) is then obtained by:
(i) forming the composition $G \circ D$,
(ii) transforming to normal form via the linear change of variables,
\begin{equation}
\left[ \begin{array}{c}
x \\ y \\ \mu
\end{array} \right] =
\frac{1}{\hat{a}_{12}^2 c^2}
\left[ \begin{array}{ccc}
1 & 0 & 0 \\
-\hat{a}_{22} & \hat{a}_{12} & \hat{b}_1 \\
0 & 0 & (1-\hat{a}_{22}) \hat{b}_1 + \hat{a}_{12} \hat{b}_2
\end{array} \right]
\left[ \begin{array}{c}
u_1 \\ w_1 \\ \eta
\end{array} \right] \;,
\label{eq:coordinateChange}
\end{equation}
where
\begin{equation}
\hat{A} = \left[ \begin{array}{cc}
\hat{a}_{11} & \hat{a}_{12} \\
\hat{a}_{21} & \hat{a}_{22}
\end{array} \right] \;, \qquad
\hat{b} = \left[ \begin{array}{c} \hat{b}_1 \\ \hat{b}_2 \end{array} \right] \;,
\end{equation}
and (iii) omitting higher order terms in $x$, $y$ and $\mu$\removableFootnote{
Moreover,
\begin{equation}
\tau = {\rm trace}(\hat{A}) \;, \qquad
\delta = \det(\hat{A}) \;, \qquad
\chi = {\rm sgn}(\hat{a}_{12} c) \;.
\end{equation}
The coordinate change cannot be performed if
(i) $c = 0$, i.e.~the square-root term is absent, or
(ii) $\hat{a}_{12} = 0$, in which case the fixed point at grazing has an eigenspace
tangent to the switching manifold.
}.
The Nordmark map ignores linear terms that are dominated by the square-root term.
For a careful consideration of these terms see \cite{MoDe01}.

\subsection{Addition of noise}
\label{sub:ADD}

The stochastic differential equation:
\begin{equation}
d\bu = \left\{ \begin{array}{lc}
f^{(L)}(\bu;\eta) \;, & u < 0 \\
f^{(R)}(\bu;\eta) \;, & u > 0 \\
\end{array} \right\} \,dt
+ \ee B(\bu;\eta) \,dW(t) \;,
\label{eq:SDE}
\end{equation}
where $W(t)$ is a vector of independent standard Brownian motions,
$B$ is matrix with smooth dependency on $\bu$ and $\eta$,
and $0 < \ee \ll 1$,
represents the system (\ref{eq:ODE})
in the presence of small-amplitude, white, Gaussian noise.
In this section we extend the construction described in the previous section
to accommodate the noise term and obtain (\ref{eq:P}).
We assume $\ee$ is small enough that the basic global structure
of the system is not destroyed by the noise and
sample paths of (\ref{eq:SDE}) likely stay near
to the associated deterministic solution, over an $O(1)$ time frame.
Theoretical results on the existence and uniqueness of solutions
to stochastic differential equations with discontinuous drift
may be found in \cite{PrSh98}, and references within.

Our strategy is to obtain stochastic versions of
the global map, $G$, (\ref{eq:G}),
and the discontinuity map, $D$, (\ref{eq:D}),
and compose them to arrive at the stochastic Nordmark map, (\ref{eq:P}).
We first discuss the stochastic global map
and consider solely the smooth left half-system of (\ref{eq:SDE}).
Issues relating to switching are discussed below.

When a sample path of the left half-system of (\ref{eq:SDE}) passes through $\Pi'$,
if the noise is not purely tangent to $\Pi'$,
with probability $1$ the sample path has multiple intersections with $\Pi'$,
albeit over an extremely short time frame.
We define the stochastic global map
by the first reintersection, $\bu_6$, of the sample path with $\Pi'$
after one large excursion.
(An alternative is to use the mean passage time \cite{Ka90}.
For stochastic Poincar\'{e} maps for systems with coloured noise, see \cite{WeKn90}.)
Multiple intersections occur sufficiently quickly
that this specification does not affect the lowest order noise term in the map \cite{FrWe84}.

For $\ee > 0$, the quantity $\bu_6$ is stochastic and
to calculate the stochastic global map
means to compute the probability density function for $\bu_6$.
We can do this to $O(\ee)$ by applying
standard asymptotic results  usually stated in the context of an exit problem.
Here we follow the sample path methodology of Freidlin and Wentzell \cite{FrWe84},
see also \cite{Sc10,GrVa99}.

For the general stochastic differential equation
\begin{equation}
d\bv = f(\bv) \,dt + \ee B(\bv) \,dW(t) \;, \qquad \bv(0) = \bv_0 \;,
\label{eq:SDEFW}
\end{equation}
where $f$ and $B$ are $C^2$,
by Theorem 2.2 of Chapter 2 of \cite{FrWe84}
we can write
\begin{equation}
\bv(t) = \bv^{(0)}(t) + \ee \bv^{(1)}(t) + o(\ee) \;,
\label{eq:xExpFW}
\end{equation}
where $\bv^{(0)}$ is the deterministic solution to (\ref{eq:SDEFW}) with $\ee=0$, 
and $\bv^{(1)}$ is the solution to the time-dependent Ornstein-Uhlenbeck process:
\begin{equation}
d\bv^{(1)} = D_\bv f(\bv^{(0)}) \bv^{(1)} \,dt + B(\bv^{(0)}) \,dW(t) \;,
\qquad \bv^{(1)}(0) = 0 \;.
\label{eq:OUFW}
\end{equation}
We then have the following lemma.

\begin{lemma}[Theorem 2.3 of Chapter 2 of Freidlin and Wentzell \cite{FrWe84}]~\\
Let $\mathcal{D}$ be a domain containing $\bv_0$
and let $T^{(0)}$ be the first passage time of $\bv^{(0)}$ from $\mathcal{D}$.
Suppose $\partial \mathcal{D}$ is differentiable at $\bv^{(0)}(T^{(0)})$,
let $p$ denote the exterior unit normal vector to $\partial \mathcal{D}$ at this point,
let $q = f \left( \bv^{(0)}(T^{(0)}) \right)$,
and assume $p^{\sf T} q \ne 0$.
If $T$ denotes the first passage time of $\bv$ from $\mathcal{D}$, then
\begin{eqnarray}
T &=& T^{(0)} - \frac{p^{\sf T} \bv^{(1)}(T^{(0)})}{p^{\sf T} q} \ee + o(\ee) \;,
\label{eq:passageTimeExpFW} \\
\bv(T) &=& \bv^{(0)}(T^{(0)}) +
\left( I - \frac{q p^{\sf T}}{p^{\sf T} q} \right) \bv^{(1)}(T^{(0)}) \ee + o(\ee) \;.
\label{eq:passagePointExpFW} 
\end{eqnarray}
\label{le:Thm23FW}
\end{lemma}


The key quantity, $\bv^{(1)}(T^{(0)})$, is stochastic, and by (\ref{eq:OUFW}),
\begin{equation}
\bv^{(1)}(t) = \int_0^t {\rm Exp} \left( \int_s^t D_\bv
f \left( \bv^{(0)}(r) \right) \,dr \right)
B \left( \bv^{(0)}(s) \right) \,dW(s) \;.
\end{equation}
Consequently $\bv^{(1)}(t)$ is a Gaussian random variable with zero mean
and covariance matrix, call it $\Omega(t)$, given by
\begin{equation}
\Omega(t) = \int_0^t H(s,t) H(s,t)^{\sf T} \,ds \;,
\label{eq:Omega}
\end{equation}
where
\begin{equation}
H(s,t) = {\rm Exp} \left( \int_s^t D_\bv f \left( \bv^{(0)}(r) \right) \,dr \right)
B \left( \bv^{(0)}(s) \right) \;.
\label{eq:H}
\end{equation}
Therefore, to $O(\ee)$, the location of first exit, $\bv(T) \in \partial \mathcal{D}$,
is Gaussian with mean, $\bv^{(0)}(T^{(0)})$, and covariance matrix
\begin{equation}
\left( I - \frac{q p^{\sf T}}{p^{\sf T} q} \right)
\Omega(T^{(0)}) \left( I - \frac{q p^{\sf T}}{p^{\sf T} q} \right)^{\sf T} \;.
\label{eq:Varw1FW}
\end{equation}

In the context of the global map, $f = f^{(L)}$, and 
$\partial \mathcal{D}$ represents $\Pi'$.
By Lemma \ref{le:Thm23FW},
\begin{equation}
\left[ \begin{array}{c}
u_6 \\ w_6
\end{array} \right] =
G(u_4,w_4;\eta)
+ \ee \xi_G + o(\ee) \;,
\label{eq:Gstoch}
\end{equation}
where $G$ is the deterministic global map, (\ref{eq:G}),
and $\xi_G$ is a zero-mean, two-dimensional Gaussian random variable.
Here $p = [0,-1,0]^{\sf T}$ and $q = [0,-\gamma_L,\zeta_L]^{\sf T} + O(|u_4,w_4,\eta|^1)$,
where $f^{(L)} \left( [0,0,0]^{\sf T} ; 0 \right) = [0,-\gamma_L,\zeta_L]^{\sf T}$.
Consequently, by (\ref{eq:Varw1FW}), the covariance matrix for $\xi_G$ is
\begin{equation}
\left[ \begin{array}{cc}
\omega_{11} & \frac{\zeta_L}{\gamma_L} \omega_{12} + \omega_{13} \\
\frac{\zeta_L}{\gamma_L} \omega_{12} + \omega_{13} &
\frac{\zeta_L^2}{\gamma_L^2} \omega_{22} + \frac{2 \zeta_L}{\gamma_L} \omega_{23} + \omega_{33}
\end{array} \right] + O(|u_4,w_4,\eta|^1) \;,
\label{eq:ThetaG}
\end{equation}
where $\omega_{ij}$ denotes the $(i,j)$-element of $\Omega(T^{(0)})$, (\ref{eq:Omega}),
and $T^{(0)}$ denotes the period of the grazing periodic orbit at $\eta = 0$\removableFootnote{
Also,
\begin{equation}
I - \frac{q p^{\sf T}}{p^{\sf T} q} =
\left[ \begin{array}{ccc}
1 & 0 & 0 \\
0 & 0 & 0 \\
0 & \frac{\zeta_L}{\gamma_L} & 1
\end{array} \right] + O(|u_4,w_4,\eta|^1) \;.
\end{equation}
We omitted the middle row and columns of the three-dimensional matrix (\ref{eq:Varw1FW})
(which are all zeros)
to obtain the two-dimensional matrix (\ref{eq:ThetaG}).
}.

A stochastic version of the discontinuity map is not as straightforward.
For instance, consider the motion from $\bu_2$ to $\bu_3$, Fig.~\ref{fig:overallSchem}.
In the absence of noise, $\bu_3$ is easily calculated using
the right half-flow.
However, with noise acting non-tangentially to $\Sigma$,
with probability $1$, a sample path from $\bu_2$ to $\bu_3$ spends a nonzero amount of time
in the left half-space.
Thus motion from $\bu_2$ to $\bu_3$ is in part governed by the left half-system, albeit weakly.
Moreover, to calculate the motion from
$\bu_1$ to $\bu_2$ and from $\bu_3$ to $\bu_4$,
we must solve the stochastic differential equation backwards in time.

Here we argue a bound for the noise in the discontinuity map.
Intuitively, the noise in the discontinuity map is greater for larger values of $u_1 > 0$.
For fixed $\ee > 0$, consider $u_1 \gg \ee$,
but recall we assume $u_1$ is small.
In this case the excursion from $\bu_1$ to $\bu_4$
is well-approximated by three
exit problems for which Lemma \ref{le:Thm23FW} is applicable,
because the weak noise does not significantly disrupt the three pieces of the excursion.
The discontinuity map acts over an $O(\sqrt{u_1})$ time-frame,
and note that a consequence of Lemma \ref{le:Thm23FW} is that for small $t$,
the standard deviations of components of the exit point are $O(\ee \sqrt{t})$.
Therefore we expect deviations in $\bu_4$ to be at most
$O \left( \ee u_1^{\frac{1}{4}} \right)$, and write
\begin{equation}
\left[ \begin{array}{c}
u_4 \\ w_4
\end{array} \right] =
D(u_1,w_1;\eta)
+ \ee \xi_D + o(\ee) \;,
\label{eq:Dstoch}
\end{equation}
where $\xi_D \le O \left( u_1^{\frac{1}{4}} \right)$.

When we compose (\ref{eq:Gstoch}) and (\ref{eq:Dstoch}),
assuming $u_1$ is small, $\xi_G$ dominates $\xi_D$ because $\xi_G = O(1)$.
As in \S\ref{sub:CONSTR},
to this composition we apply the change of variables (\ref{eq:coordinateChange}),
and omit higher order terms, producing (\ref{eq:P}).
Moreover, the noise term in (\ref{eq:P}) is
$\xi = \frac{1}{\hat{a}_{12}^2 c^2} \left[ \begin{array}{cc}
1 & 0 \\
-\hat{a}_{22} & \hat{a}_{12}
\end{array} \right] \xi_G$,
and therefore the covariance matrix for $\xi$ is given by
\begin{equation}
\Theta =
\frac{1}{\hat{a}_{12}^4 c^4}
\left[ \begin{array}{cc}
1 & 0 \\
-\hat{a}_{22} & \hat{a}_{12}
\end{array} \right]
\left[ \begin{array}{cc}
\omega_{11} & \frac{\zeta_L}{\gamma_L} \omega_{12} + \omega_{13} \\
\frac{\zeta_L}{\gamma_L} \omega_{12} + \omega_{13} &
\frac{\zeta_L^2}{\gamma_L^2} \omega_{22} + \frac{2 \zeta_L}{\gamma_L} \omega_{23} + \omega_{33}
\end{array} \right]
\left[ \begin{array}{cc}
1 & -\hat{a}_{22} \\
0 & \hat{a}_{12}
\end{array} \right] \;.
\label{eq:Theta}
\end{equation}

\subsection{A linear oscillator with compliant impacts}
\label{sub:IMPACT}

\begin{figure}[b!]
\begin{center}
\setlength{\unitlength}{1cm}
\begin{picture}(10.6,4.1)
\put(0,.1){\includegraphics[height=4cm]{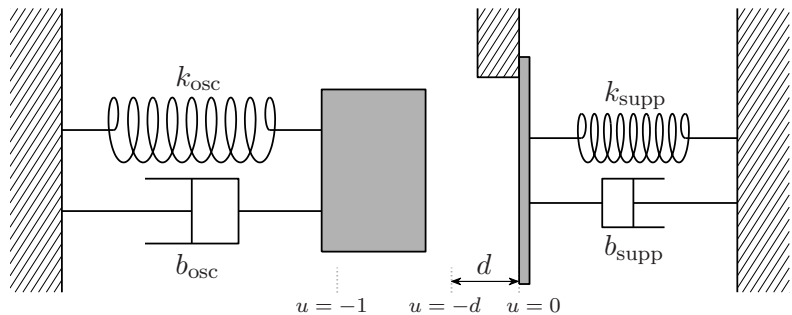}}	
\put(6.8,0){\scriptsize $u=0$}
\put(5.5,0){\scriptsize $u=-d$}
\put(4,0){\scriptsize $u=-1$}
\put(6.4,.5){\small $d$}
\put(2.4,3){\small $k_{\rm osc}$}
\put(2.4,.5){\small $b_{\rm osc}$}
\put(8.1,2.8){\small $k_{\rm supp}$}
\put(8.1,.75){\small $b_{\rm supp}$}
\end{picture}
\caption{
A schematic of the forced vibro-impacting system modelled by (\ref{eq:impactOsc}).
The equilibrium positions of the mass and support are $u=-1$ and $u=-d$ respectively.
\label{fig:vibroImpactCompliant}
}
\end{center}
\end{figure}

As an example, we consider a harmonically forced linear oscillator experiencing compliant impacts
with a massless, prestressed support.
The system is depicted in Fig.~\ref{fig:vibroImpactCompliant}
and we model it with the non-dimensionalized equation of motion\removableFootnote{
Can we call this a Langevin equation?
Here $W$ is scalar, whereas everywhere else in the paper $W$ is a vector.
Would it be better to use a new symbol to denote $\frac{dW}{dt}$,
say $\nu$ or $\varphi$, or $\eta$ but then replace the $\eta$ that is used elsewhere?
},
\begin{equation}
\ddot{u} = \left\{ \begin{array}{lc}
-k_{\rm osc}(u+1) - b_{\rm osc} \dot{u} \;, & u < 0 \\
-k_{\rm osc}(u+1) - (b_{\rm osc} + b_{\rm supp}) \dot{u}
-k_{\rm supp}(u+d) \;, & u > 0
\end{array} \right\} + F \cos(t) + \ee \nu(t) \;,
\label{eq:impactOsc}
\end{equation}
where $\nu(t)$ is a standard Gaussian white noise process.
Here $u$ denotes the position of the mass.
A rigid stop placed at $u=0$
prevents the support moving left of this value and
prestresses the support by a distance $d>0$.
We assume that while the mass is not in contact with the support,
the support is at rest at $u=0$.
Whenever the mass arrives at $u=0$,
the mass and support then move together for $u>0$ until
a return to $u=0$ at which point the stop catches the support
and the mass returns to non-impacting motion\removableFootnote{
We make several more assumptions on the system
but omit them for brevity.
These include assuming the motion is one-dimensional.
We assume the stop does not directly interfere with the motion of the mass.
We ignore energy loss at impacts.
Finally note that it is reasonable to assume that the mass does not detach
from the support at any $u>0$ when $d>0$ and impact velocities are low,
because, it may be shown that if the ratio of impact velocity to $d$
is sufficiently small then the contact force between the mass and the support is never zero.
}.
The setup is identical to a scenario given in \cite{MaIn08}.
Experiments involving compliant impacts are described in \cite{InPa08b,InPa06}.
We incorporate randomness in the model
by adding white-noise to the forcing term.

\begin{figure}[b!]
\begin{center}
\setlength{\unitlength}{1cm}
\begin{picture}(14.9,11.85)
\put(.2,6.6){\includegraphics[height=5.25cm]{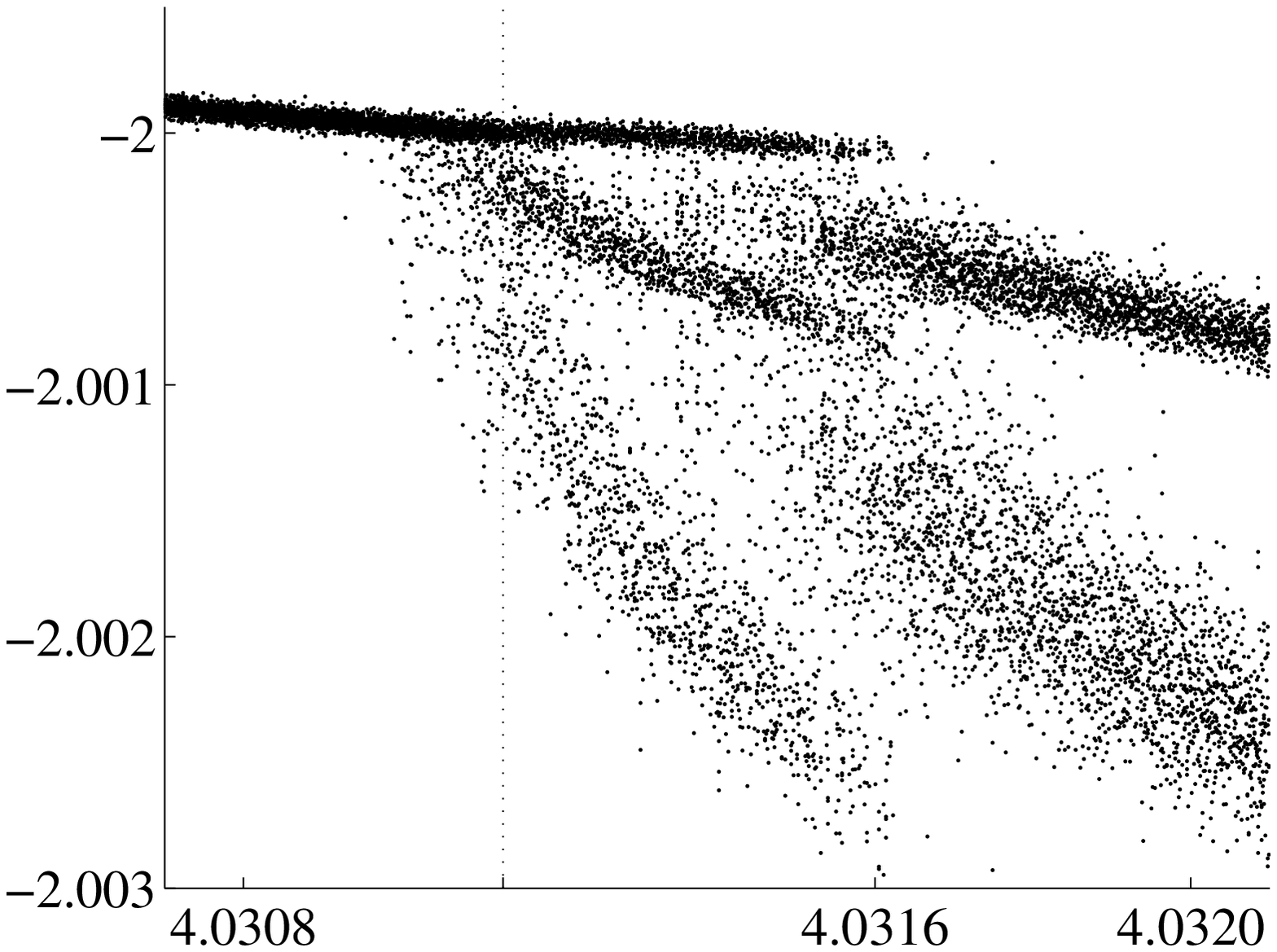}}
\put(7.9,6.6){\includegraphics[height=5.25cm]{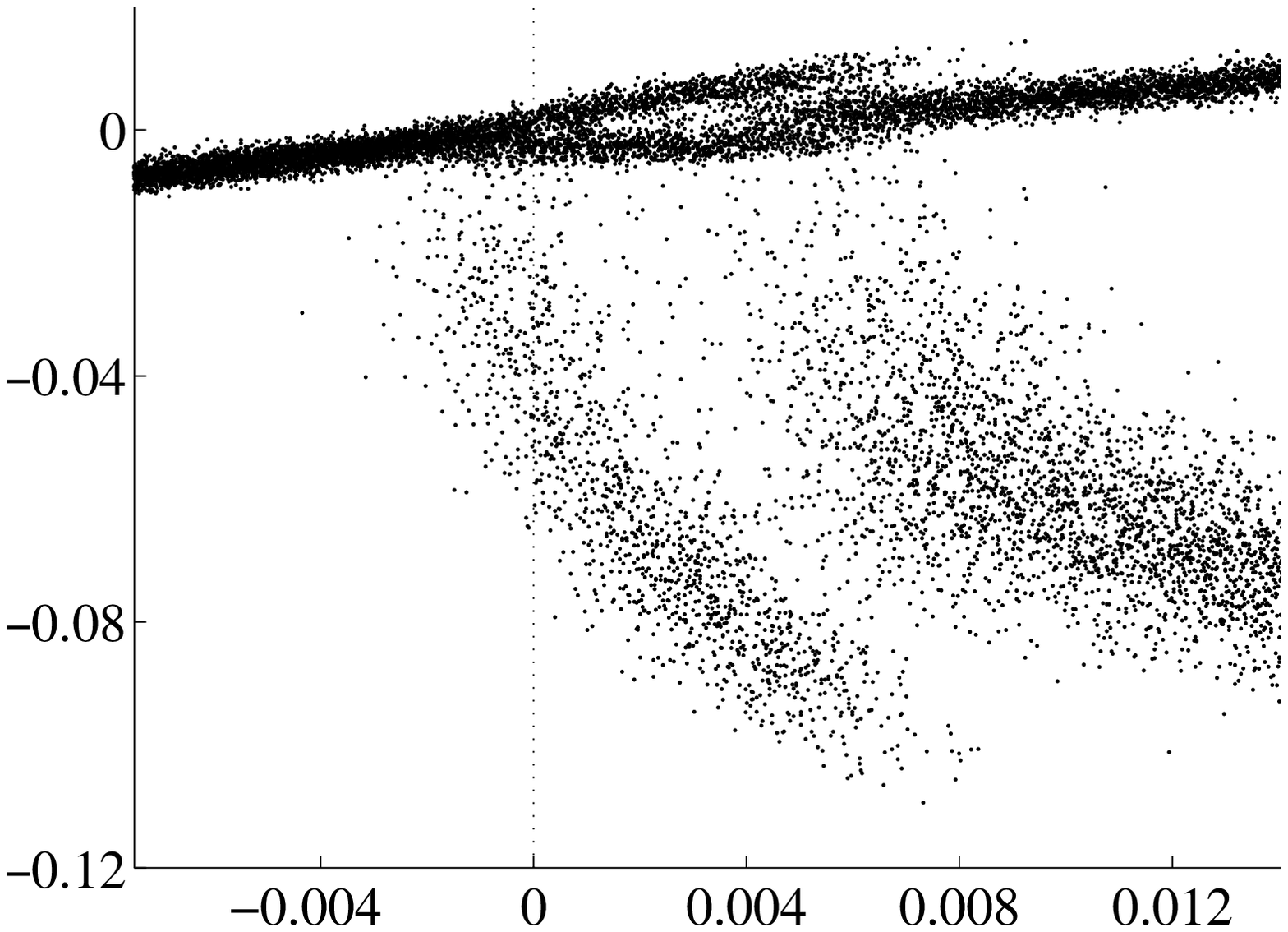}}
\put(.2,.3){\includegraphics[height=5.25cm]{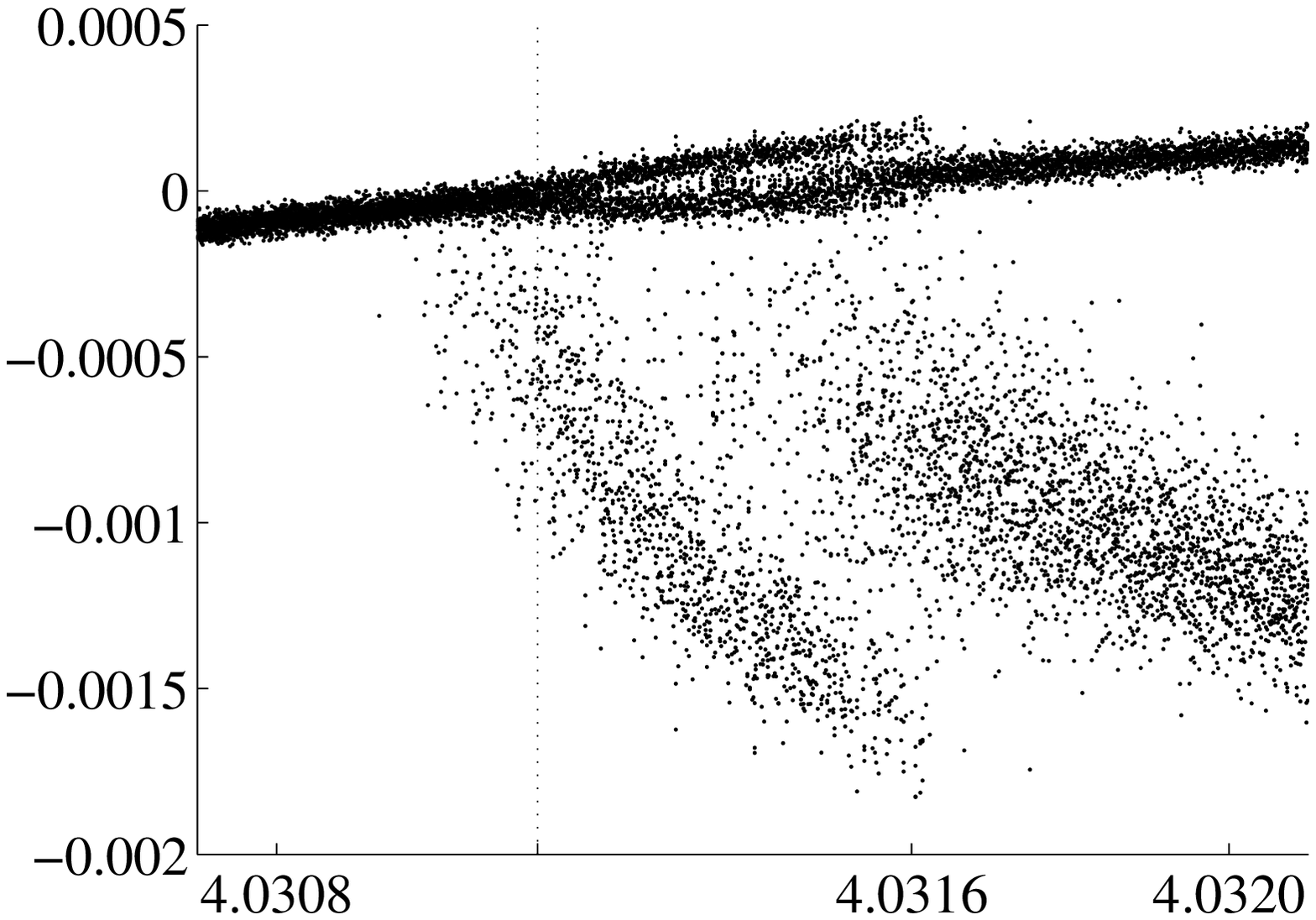}}
\put(7.9,.3){\includegraphics[height=5.25cm]{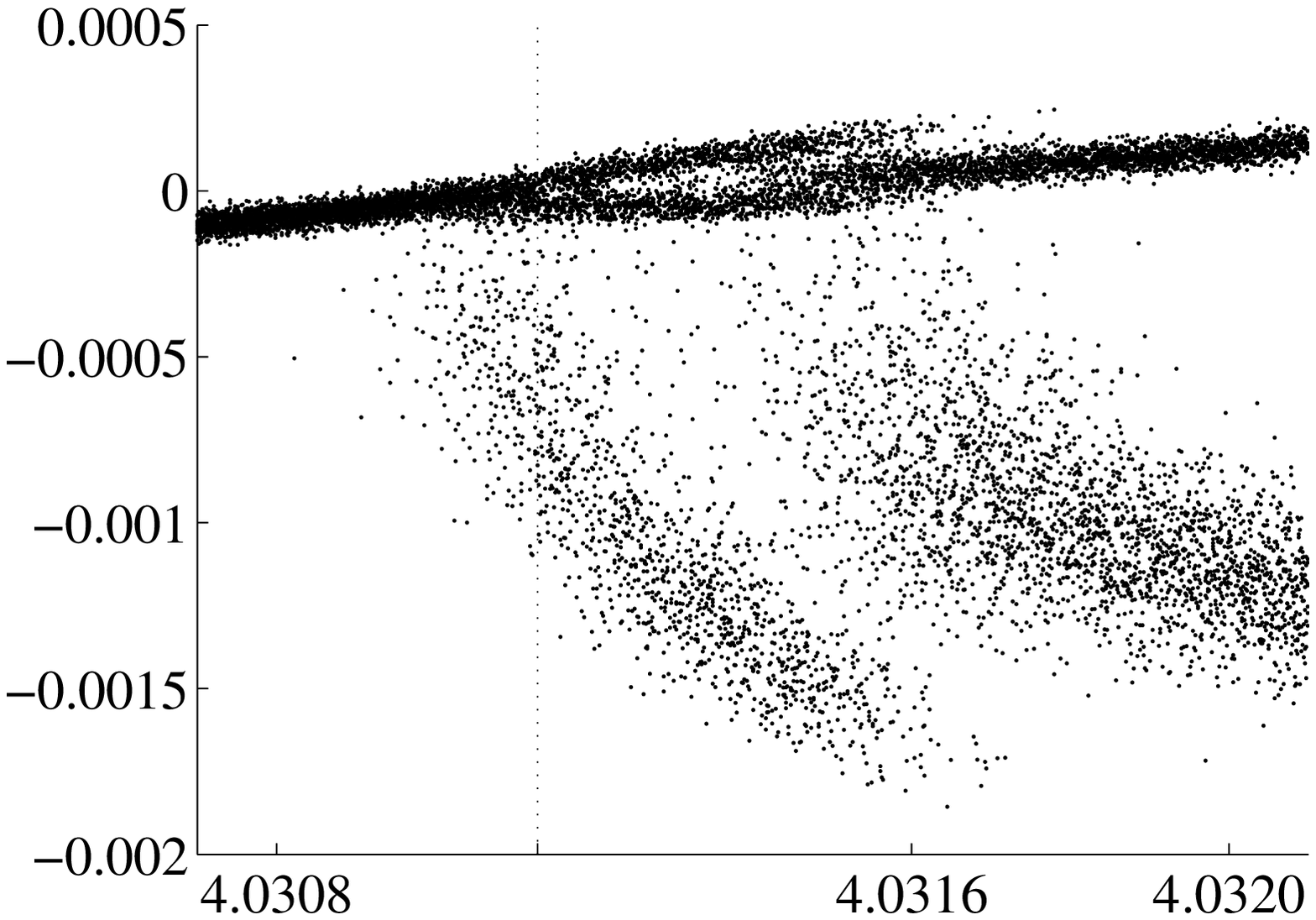}}
\put(4,6.3){$F$}
\put(0,9.1){$u_0$}
\put(2.8,6.63){\scriptsize $F_{\rm graz}$}
\put(11.7,6.3){$\mu$}
\put(7.7,9.1){$x$}
\put(4,0){$F$}
\put(0,2.8){$u_1$}
\put(2.8,.33){\scriptsize $F_{\rm graz}$}
\put(11.7,0){$F$}
\put(7.7,2.8){$u_1$}
\put(10.5,.33){\scriptsize $F_{\rm graz}$}
\put(2,11.75){\large \sf \bfseries A}
\put(9.7,11.75){\large \sf \bfseries B}
\put(2,5.45){\large \sf \bfseries C}
\put(9.7,5.45){\large \sf \bfseries D}
\end{picture}
\caption{
Panels A and C show numerically computed bifurcation diagrams
of the impact oscillator model (\ref{eq:impactOsc}) with
$(k_{\rm osc},b_{\rm osc},k_{\rm supp},b_{\rm supp},d) = (5,0.5,10,0,0.1)$
and $\ee = 5 \times 10^{-5}$.
In panel A, the values are intersections of long-time dynamics
with $\Pi = \{ (u,\dot{u}) ~|~ \dot{u}=0, u \approx -2 \}$.
Panel C shows intersections with $\Pi'$.
For $u_1 > 0$ these intersections are virtual 
and were obtained by appropriately
solving the non-impacting equation of motion for $u > 0$.
Panel B is a bifurcation diagram of the stochastic Nordmark map, (\ref{eq:P}),
using parameter values corresponding to the system in panel A
as calculated by the formulas given in the text.
Specifically,
$(\tau,\delta) \approx (0.07264,0.04321)$, $\chi = 1$, and
$(\theta_{11},\theta_{12},\theta_{22}) \approx (662.6,-7.450,28.29)$,
where $\theta_{ij}$ denotes the $(i,j)$-element of $\Theta$.
Panel D shows the points of panel B
under the inverse of the coordinate change (\ref{eq:coordinateChange}).
\label{fig:checkNormalSoft}
}
\end{center}
\end{figure}

To simplify our discussion we assume,
$\frac{b_{\rm osc}^2}{4} < k_{\rm osc} < 1$,
so that, in particular, the system is sub-resonant\removableFootnote{
Resonance occurs when $k_{\rm osc} = 1$.
The Floquet multipliers associated with the non-impacting
periodic solution are complex-valued when $k_{\rm osc} > \frac{b_{\rm osc}^2}{4}$.
Also $\frac{\pi}{2} < t_{\rm graz} < \pi$
if and only if $k_{\rm osc} < 1$.
}.
We treat the forcing amplitude, $F$, as the primary bifurcation parameter.
When $\ee = 0$, for small $F$ the system settles to a stable non-impacting
periodic orbit of period $2 \pi$.
The periodic orbit undergoes a grazing bifurcation at
\begin{equation}
F_{\rm graz} = \sqrt{b_{\rm osc}^2 + (1-k_{\rm osc})^2} \;.
\label{eq:FGraz}
\end{equation}
A typical bifurcation diagram of (\ref{eq:impactOsc})
with $\ee > 0$ is shown in Fig.~\ref{fig:checkNormalSoft}-A.
Fig.~\ref{fig:checkNormalSoft}-B is a bifurcation diagram
of the stochastic Nordmark map, (\ref{eq:P}),
with parameter values corresponding to the impact oscillator for panel A.
As detailed in Appendix \ref{sec:APP},
the values $\tau$, $\delta$ and $\chi$,
are given by the deterministic results of \S\ref{sub:CONSTR},
and $\Theta$ is given by (\ref{eq:Theta}).
As expected, the two bifurcation diagrams exhibit the same qualitative structure.
To quantitatively compare (\ref{eq:impactOsc}) and (\ref{eq:P}),
for each we have computed the
corresponding bifurcation diagrams on $\Pi'$, panels C and D\removableFootnote{
Previously I had compared panel A with the data of panel D
under the inverse image of deterministic Poincar\'{e} map
that takes points from $\Pi$ to $\Pi'$.
Such a comparison is not valid because under the inverse image of this map
the data experiences an expansion due to the contracting nature of the map
that should be curbed the noise.
}.
As expected these bifurcation diagrams are practically indistinguishable,
justifying our use of the map (\ref{eq:P}).
A more precise comparison is beyond the scope of this paper.

\section{Approximating invariant densities by a linear analysis}
\label{sec:LINEAR}
\setcounter{equation}{0}

In this section we consider the case where the invariant density can be approximated
by a combination of Gaussians and calculate the related covariance matrices.
The calculations reveal instances of larger variability
and situations where the Gaussian approximation may break down.

We consider the case that (\ref{eq:P}) with $\ee = 0$
has an attracting periodic solution, $\{ (x_i,y_i) \}_{i = 0,\ldots,n-1}$,
with $x_i \ne 0$, for all $i$.
Then with small $\ee > 0$, numerical simulations indicate that (\ref{eq:P}) has a
stable invariant density
with high peaks near each $(x_i,y_i)$.
For sufficiently small $\ee > 0$,
the bulk of this invariant density is distant from the switching manifold.
Consequently in this case the nonsmooth effect of the switching manifold is not seen.
About each $(x_i,y_i)$ the invariant density has a roughly Gaussian shape
that may be computed by linearizing the map about $(x_i,y_i)$.
Therefore the invariant density of (\ref{eq:P})
is well-approximated by the sum of $n$ Gaussian densities centred at each $(x_i,y_i)$,
and scaled by $\frac{1}{n}$ such that the density is normalized
(see Fig.~\ref{fig:invDensity_a} for an example).
In \S\ref{sub:MPS} we derive an appropriate $n^{\rm th}$ iterate of (\ref{eq:P})
and solve for $(x_i,y_i)$ when $\ee = 0$.
In \S\ref{sub:GAUSS} we calculate the covariance matrices of the Gaussians.
Further calculations and approximations are described in
\S\ref{sub:APPROX} and \S\ref{sub:THETA}.


\subsection{Maximal periodic solutions}
\label{sub:MPS}

Stability of a periodic solution, $\{ (x_i,y_i) \}_{i = 0,\ldots,n-1}$,
of (\ref{eq:P}) with $\ee = 0$
may be determined in a standard manner from the multipliers of the matrix
formed from the product of Jacobians evaluated at each $(x_i,y_i)$.
For each point of the solution for which $x_i > 0$,
the product accumulates a term of order $x_i^{-\frac{1}{2}}$,
which is strongly destabilizing since we assume the $x_i$ are small.
For this reason we only consider periodic solutions
for which exactly one iterate satisfies $x_i > 0$.
We expect no other periodic solutions to have a significant effect
on long-term dynamics near the grazing bifurcation.
Indeed this is consistent with numerical results.
This observation has been noted previously and 
such solutions are known as {\em maximal periodic solutions} \cite{No01,ChOt94}.

Suppose that the next $n-1$ iterates of (\ref{eq:P})
from a point $(x_0,y_0)$ with $x_0 \ge 0$ satisfy $x_i \le 0$.
Then the $n^{\rm th}$-iterate is given by
\begin{eqnarray}
\left[ \begin{array}{c} x_n \\ y_n \end{array} \right]
&=& A^n
\left[ \begin{array}{c} x_0 \\ y_0 - \chi \sqrt{x_0} \end{array} \right]
+ (I + A + \cdots + A^{n-1})
\left[ \begin{array}{c} 0 \\ 1 \end{array} \right] \mu \nonumber \\
&&+~\ee
\left( \xi_{n-1} + A \xi_{n-2} + \cdots +
A^{n-1} \xi_0 \right) \;,
\label{eq:nthIterateStoch}
\end{eqnarray}
where the $n$ random vectors, $\xi_i$, are independent.
The covariance matrix of the sum of $n$ independent Gaussian random variables
the sum of the covariance matrices of each\removableFootnote{
This can be demonstrated with moment generating functions.
},
therefore
\begin{equation}
\xi_{n-1} + A \xi_{n-2} + \cdots + A^{n-1} \xi_0 \sim \xi^{(n)} \;,
\end{equation}
where $\xi^{(n)}$ is a zero-mean Gaussian random variable with covariance matrix
\begin{equation}
\Theta^{(n)} = \sum_{i=0}^{n-1} A^i \Theta \left( A^{\sf T} \right)^i \;.
\label{eq:nthIterateCov2}
\end{equation}
Hence we can write (\ref{eq:nthIterateStoch}) as
\begin{equation}
\left[ \begin{array}{c} x_n \\ y_n \end{array} \right]
= A^n \left[ \begin{array}{c} x_0 \\ y_0 - \chi \sqrt{x_0} \end{array} \right]
+ b^{(n)} \mu + \ee \xi^{(n)} \;,
\label{eq:nthIterateStoch2}
\end{equation}
where
\begin{equation}
b^{(n)} = (I + A + \cdots + A^{n-1})
\left[ \begin{array}{c} 0 \\ 1 \end{array} \right] \;.
\end{equation}
Fixed points, $(x^{*(n)},y^{*(n)})$, of (\ref{eq:nthIterateStoch2}) with $\ee = 0$, satisfy
\begin{eqnarray}
y^{*(n)} &=& \frac{(a_{12} b_2 - a_{22} b_1) \chi \sqrt{x^{*(n)}} +
\left( (1-a_{11}) b_2 + a_{21} b_1 \right) x^{*(n)}}
{(1-a_{22}) b_1 + a_{12} b_2} \;, \label{eq:ynxn} \\
\mu &=& \frac{a_2 \chi \sqrt{x^{*(n)}} +
\left( (1-a_{11}) (1-a_{22}) - a_{12} a_{21} \right) x^{*(n)}}
{(1-a_{22}) b_1 + a_{12} b_2} \;, \label{eq:muxn}
\end{eqnarray}
where
\begin{equation}
A^n = \left[ \begin{array}{cc}
a_{11} & a_{12} \\ a_{21} & a_{22}
\end{array} \right] \;, \qquad
b^{(n)} = \left[ \begin{array}{c}
b_1 \\ b_2
\end{array} \right] \;.
\label{eq:Jnfnbn}
\end{equation}
The stability multipliers associated with $(x^{*(n)},y^{*(n)})$ are
the eigenvalues of the matrix
\begin{equation}
K(x^{*(n)}) = A^n \left[ \begin{array}{cc}
1 & 0 \\ \frac{-\chi}{2 \sqrt{x^{*(n)}}} & 1
\end{array} \right] =
\left[ \begin{array}{cc}
a_{11} - \frac{a_{12} \chi}{2 \sqrt{x^{*(n)}}} & a_{12} \\ a_{21}
- \frac{a_{22} \chi}{2 \sqrt{x^{*(n)}}} & a_{22}
\end{array} \right] =
\left[ \begin{array}{cc}
k_{11} & k_{12} \\
k_{21} & k_{22}
\end{array} \right] \;.
\label{eq:K}
\end{equation}

In summary, the point of a maximal periodic solution with $x > 0$,
$(x^{*(n)},y^{*(n)})$, is given implicitly in terms of $\mu$
by (\ref{eq:ynxn}) and (\ref{eq:muxn}).
The remaining points are easily obtained by iterating (\ref{eq:P}) from this point.
The periodic solution is stable when both eigenvalues of $K(x^{*(n)})$, (\ref{eq:K}),
lie inside the unit circle.
The blue curves in
Figs.~\ref{fig:stochBifDiag_c} and \ref{fig:stochBifDiag_b}
denote the $x$-values of stable, admissible, maximal periodic solutions
(up to period $5$) and were computed by evaluating the above expressions.

\subsection{Calculation of Gaussian invariant densities}
\label{sub:GAUSS}

Here we determine the covariance matrices for the Gaussian approximations.
We first compute the covariance matrix for the Gaussian centred at $(x^{*(n)},y^{*(n)})$,
then give an iterative formula for the remaining covariance matrices.

The linearization of the $n^{\rm th}$-iterate, (\ref{eq:nthIterateStoch2}),
about the fixed point, $(x^{*(n)},y^{*(n)})$, is
\begin{equation}
\bz_{j+1} = K \bz_j + \ee \xi^{(n)} \;,
\label{eq:nthLinear}
\end{equation}
where $K = K(x^{*(n)}(\mu))$, (\ref{eq:K}), and\removableFootnote{
Consider the stochastic map
\begin{equation}
x_{i+1} = S(x_i) + \xi_i \;,
\label{eq:nlsm}
\end{equation}
where $x_i \in \mathbb{R}^N$,
$S$ is a non-stochastic function
and the $\xi_i$ are independent random vectors having a common PDF $g(x)$.
(For linear maps for which the functional part $S$ is stochastic, see \cite{Br86,Ve79}.)
It is instructive to calculate how (\ref{eq:nlsm}) transforms densities.
Let $f_i(x)$ denote the density after $i$ iterations of (\ref{eq:nlsm})
from some initial PDF, $f_0(x)$.
By conditioning over all possible previous points we obtain the formula
\begin{equation}
f_{i+1}(x) = \int_{\mathbb{R}^N} f_i(y) g(x - S(y)) \,dy \;,
\label{eq:densityIteration}
\end{equation}
(see \cite{LaMa94} for a formal derivation).
If (\ref{eq:nlsm}) is deterministic, say $g(x) = \delta(x)$,
and $S$ is one-to-one, then (\ref{eq:densityIteration}) simplifies to
\begin{equation}
f_{i+1}(x) = f_i(S^{-1}(x)) | D_x S^{-1}(x) | \;,
\label{eq:FrobeniusPerron}
\end{equation}
which is the Frobenius-Perron operator \cite{LaMa94,DiZh09}
that describes transformations of densities under $S$.
If $\xi_i$ in (\ref{eq:nlsm}) is replaced with $\varepsilon \xi_i$,
then the operator (\ref{eq:densityIteration}) limits to (\ref{eq:FrobeniusPerron})
as $\varepsilon \to 0$ \cite{LaMa94}.
Of special interest are densities $f$ that are invariant under (\ref{eq:nlsm}).
If $S$ is a contracting map,
then (\ref{eq:nlsm}) has a unique stable invariant density \cite{LaMa94}
(Uniqueness is given by Theorem 10.5.2;
existence and stability are given by Example 10.5.2.).
}
\begin{equation}
\bz = \left[ \begin{array}{c}
x - x^{*(n)}(\mu) \\
y - y^{*(n)}(\mu)
\end{array} \right] \;,
\end{equation}
denotes the displacement from the fixed point.
Since (\ref{eq:nthLinear}) is linear and $\xi^{(n)}$
has a Gaussian distribution (with zero mean and covariance matrix $\Theta^{(n)}$,
(\ref{eq:nthIterateCov2})),
(\ref{eq:nthLinear}) maps one Gaussian probability density function
of $\bz$-values, into another \cite{So02,LaMa94}.
Specifically, if the $j^{\rm th}$ density under (\ref{eq:nthLinear})
is Gaussian with zero mean and covariance matrix $\ee^2 \Lambda_j$, then
\begin{equation}
\Lambda_{j+1} = K \Lambda_j K^{\sf T} + \Theta^{(n)} \;.
\label{eq:densityMap}
\end{equation}
If both eigenvalues of $K$ lie inside the unit circle,
then the fixed point $(x^{*(n)},y^{*(n)})$ of (\ref{eq:nthIterateStoch2}) is attracting
and (\ref{eq:nthLinear}) converges to a unique, zero-mean, Gaussian, invariant density
with covariance matrix, $\ee^2 \Lambda$, satisfying
\begin{equation}
\Lambda = K \Lambda K^{\sf T} + \Theta^{(n)} \;.
\label{eq:Lyapunov}
\end{equation}

Equation (\ref{eq:Lyapunov}) is a {\em Lyapunov equation}
that cannot be solved for $\Lambda$ by
elementary matrix algebra operations \cite{So02}.
However we can break (\ref{eq:Lyapunov}) into
its constituent scalar equations\removableFootnote{
$N$-step scalar autoregression models
(used to model time series data)
may be written as $N$-dimensional linear stochastic maps \cite{So02,Ha94}.
In this case the matrix $K$ is a companion matrix
and the so-called {\em Yule-Walker equations}
are used to compute the covariance matrix of the invariant density.
An $N$-dimensional Lyapunov equation is broken apart into $N^2$ linear scalar equations.
Since covariance matrices are symmetric,
$\frac{N(N-1)}{2}$ of these equations will be redundant 
so the elements of $\Lambda$ may be written
solution to an $\frac{N(N+1)}{2}$-dimensional matrix equation.
}.
We write (\ref{eq:Lyapunov}) as
\begin{equation}
\left[ \begin{array}{cc}
\lambda_{11} & \lambda_{12} \\
\lambda_{12} & \lambda_{22}
\end{array} \right] =
\left[ \begin{array}{cc}
k_{11} & k_{12} \\
k_{21} & k_{22}
\end{array} \right]
\left[ \begin{array}{cc}
\lambda_{11} & \lambda_{12} \\
\lambda_{12} & \lambda_{22}
\end{array} \right]
\left[ \begin{array}{cc}
k_{11} & k_{21} \\
k_{12} & k_{22}
\end{array} \right] +
\left[ \begin{array}{cc}
\theta^{(n)}_{11} & \theta^{(n)}_{12} \\
\theta^{(n)}_{12} & \theta^{(n)}_{22}
\end{array} \right] \;,
\label{eq:Lyapunov2d}
\end{equation}
which may be broken up and rewritten as
\begin{equation}
\left[ \begin{array}{c}
\lambda_{11} \\ \lambda_{12} \\ \lambda_{22}
\end{array} \right] = M
\left[ \begin{array}{c}
\lambda_{11} \\ \lambda_{12} \\ \lambda_{22}
\end{array} \right] +
\left[ \begin{array}{c}
\theta^{(n)}_{11} \\ \theta^{(n)}_{12} \\ \theta^{(n)}_{22}
\end{array} \right] \;,
\end{equation}
where
\begin{equation}
M = \left[ \begin{array}{ccc}
k_{11}^2 & 2 k_{11} k_{12} & k_{12}^2 \\
k_{11} k_{21} & k_{11} k_{22} + k_{12} k_{21} & k_{12} k_{22} \\
k_{21}^2 & 2 k_{21} k_{22} & k_{22}^2
\end{array} \right] \;.
\label{eq:M}
\end{equation}
Then
\begin{equation}
\left[ \begin{array}{c}
\lambda_{11} \\ \lambda_{12} \\ \lambda_{22}
\end{array} \right] = (I-M)^{-1}
\left[ \begin{array}{c}
\theta^{(n)}_{11} \\ \theta^{(n)}_{12} \\ \theta^{(n)}_{22}
\end{array} \right] \;.
\label{eq:lambdas}
\end{equation}
A direct calculation reveals\removableFootnote{
Also
\begin{equation}
{\rm adj}(I-M) = \left[ \begin{array}{ccc}
1 - k_{11} k_{22} - k_{12} k_{21} - k_{22}^2 (1 - \delta^n) &
2 k_{12} \left( k_{11} - k_{22} \delta^n \right) &
k_{12}^2 (1 + \delta^n) \\
k_{21} \left( k_{11} - k_{22} \delta^n \right) &
(1 - k_{11}^2)(1 - k_{22}^2) - k_{12}^2 k_{21}^2 &
k_{12} \left( k_{22} - k_{11} \delta^n \right) \\
k_{21}^2 (1 + \delta^n) &
2 k_{21} \left( k_{22} - k_{11} \delta^n \right) &
1 - k_{11} k_{22} - k_{12} k_{21} - k_{11}^2 (1 - \delta^n)
\end{array} \right] \;,
\end{equation}
where we have substituted
$k_{11} k_{22} - k_{12} k_{21} = \det(K) = \det(A^n) = \det(A)^n = \delta^n$.
}
\begin{equation}
\det(I-M) =
\big( \det(K) - {\rm trace}(K) + 1 \big)
\big( \det(K) + {\rm trace}(K) + 1 \big)
\big( 1 - \det(K) \big) \;,
\label{eq:detImM}
\end{equation}
thus $\det(I-M) = 0$ if and only if $K$ has an eigenvalue on the unit circle.
Moreover, if under parameter change,
an eigenvalue of $K$ approaches the unit circle,
$\Lambda$ tends to infinity.

$\Lambda$, given by explicitly by
(\ref{eq:lambdas}) with (\ref{eq:M}) and (\ref{eq:K}),
is the covariance matrix for the Gaussian approximation to
the invariant density about $(x^{*(n)},y^{*(n)})$.
The covariance matrix relating to the $i^{\rm th}$ iterate of
$(x^{*(n)},y^{*(n)})$, call it $\Lambda^{(i)}$,
may be calculated via the identity
\begin{equation}
\Lambda^{(i+1)} = A \Lambda^{(i)} A^{\sf T} + \Theta \;,
\label{eq:Lambdai}
\end{equation}
(where $\Lambda^{(n)} = \Lambda^{(0)} = \Lambda$),
which is valid for $i = 1, \ldots, n-1$.
Let $\lambda_{11}^{(i)}$ denote the $(1,1)$-element of $\Lambda^{(i)}$.
The red curves in Figs.~\ref{fig:stochBifDiag_c} and \ref{fig:stochBifDiag_b}
denote a distance $\ee \lambda_{11}^{(i)}$ from the $i^{\rm th}$ iterates
of the periodic solution (the blue curves)
and correspond to one standard deviation from the mean.
The $\lambda_{11}^{(i)}$ may differ significantly for different $i$
as in Fig.~\ref{fig:checkNormalSoft}.

\subsection{An approximation to the covariance matrices}
\label{sub:APPROX}

Recall, the stability multipliers of a maximal periodic solution
with a point, $(x^{*(n)},y^{*(n)})$, in the right half-plane
are the eigenvalues of $K$, (\ref{eq:K}).
Near grazing, $x^{*(n)}$ is small,
thus if this periodic solution is attracting,
the elements, $a_{12}$ and $a_{22}$, of $A^n$, must also be small.
In view of (\ref{eq:K}), it is reasonable to suppose
$|a_{ij}| \le O \left( \sqrt{x^{*(n)}} \right)$, for each $a_{ij}$.
Then from (\ref{eq:M}),
\begin{equation}
M = \left[ \begin{array}{ccc}
\frac{a_{12}^2}{4 x^{*(n)}} & 0 & 0 \\
\frac{a_{12} a_{22}}{4 x^{*(n)}} & 0 & 0 \\
\frac{a_{22}^2}{4 x^{*(n)}} & 0 & 0
\end{array} \right] +
O \left( \sqrt{x^{*(n)}} \right) \;,
\end{equation}
and thus
\begin{equation}
(I-M)^{-1} = I + \frac{1}{1 - \frac{a_{12}^2}{4 x^{*(n)}}}
\left[ \begin{array}{ccc}
\frac{a_{12}^2}{4 x^{*(n)}} & 0 & 0 \\
\frac{a_{12} a_{22}}{4 x^{*(n)}} & 0 & 0 \\
\frac{a_{22}^2}{4 x^{*(n)}} & 0 & 0
\end{array} \right] +
O \left( \sqrt{x^{*(n)}} \right) \;.
\end{equation}
Then by (\ref{eq:lambdas}),
\begin{equation}
\Lambda = \Theta^{(n)} +
\frac{\theta^{(n)}_{11}}{4 x^{*(n)} - a_{12}^2}
\left[ \begin{array}{cc}
a_{12}^2 & a_{12} a_{22} \\
a_{12} a_{22} & a_{22}^2
\end{array} \right] +
O \left( \sqrt{x^{*(n)}} \right) \;.
\label{eq:LambdaApprox}
\end{equation}
Equation (\ref{eq:LambdaApprox})
indicates the deviation of $\Lambda$
(the covariance matrix for the Gaussian invariant density about $(x^{*(n)},y^{*(n)})$)
from $\Theta^{(n)}$, (\ref{eq:nthIterateCov2}).
For large $n$ we often have $|a_{ij}| \ll \sqrt{x^{*(n)}}$,
in which case $\Lambda \approx \Theta^{(n)}$.
Alternatively, if the periodic solution is weakly attracting,
then by (\ref{eq:detImM}), the denominator, $4 x^{*(n)} - a_{12}^2$, is small.
By (\ref{eq:LambdaApprox}), in this case
$\Lambda$ deviates significantly from $\Theta^{(n)}$,
as seen in Figs.~\ref{fig:stochBifDiag_c} and \ref{fig:stochBifDiag_b}
where the red curves are relatively distant from the associated blue curves.


\subsection{Effects of different $\Theta$}
\label{sub:THETA}

Both eigenvalues of $A$ lie inside the unit circle
(assuming the grazing periodic orbit is attracting)
and therefore the infinite series,
\begin{equation}
\Theta^{(\infty)} \equiv \sum_{i=0}^\infty A^i \Theta \left( A^{\sf T} \right)^i \;,
\label{eq:infIterateCov}
\end{equation}
converges.
Thus for large $n$, by (\ref{eq:nthIterateCov2}),
$\Theta^{(\infty)}$ may provide a good approximation to $\Theta^{(n)}$.
Here we derive a simple expression for $\Theta^{(\infty)}$.
Equation (\ref{eq:infIterateCov}) is equivalent to the Lyapunov equation,
\begin{equation}
\Theta = \Theta^{(\infty)} - A \Theta^{(\infty)} A^{\sf T} \;.
\end{equation}
By expanding the individual elements of this equation
and solving we obtain
\begin{equation}
\left[ \begin{array}{c}
\theta^{(\infty)}_{11} \\ \theta^{(\infty)}_{12} \\ \theta^{(\infty)}_{22}
\end{array} \right] =
\frac{1}{\Delta}
\left[ \begin{array}{ccc}
1+\delta & 2 \tau & 1+\delta \\
-\tau \delta & 1-\tau^2-\delta^2 & -\tau \delta \\
\delta^2+\delta^3 & 2 \tau \delta^2 & 1+\delta-\tau^2+\tau^2 \delta
\end{array} \right]
\left[ \begin{array}{c}
\theta_{11} \\ \theta_{12} \\ \theta_{22}
\end{array} \right] \;,
\label{eq:infthetas2}
\end{equation}
where
\begin{equation}
\Delta = (\delta-\tau+1)(\delta+\tau+1)(1-\delta) > 0 \;.
\end{equation}

We can use (\ref{eq:infthetas2}) to investigate the effects of different $\Theta$.
In order to keep the noise amplitude controlled by $\ee$,
here we fix ${\rm trace}(\Theta) = 2$.
Substituting $\theta_{11} = 2 - \theta_{22}$ into (\ref{eq:infthetas2}) yields
\begin{equation}
\left[ \begin{array}{c}
\theta^{(\infty)}_{11} \\ \theta^{(\infty)}_{12} \\ \theta^{(\infty)}_{22}
\end{array} \right] =
\left[ \begin{array}{c}
\frac{2}{\Delta} \left( 1 + \delta + \tau \theta_{12} \right) \\
\frac{1}{\Delta} \left( -2 \tau \delta + (1-\tau^2-\delta^2) \theta_{12} \right) \\
\frac{2 \delta^2}{\Delta} \left( 1 + \delta + \tau \theta_{12} \right) + \theta_{22}
\end{array} \right] \;,
\label{eq:infthetas3}
\end{equation}
where $\theta_{12}$ and $\theta_{22}$ are constrained by
\begin{equation}
\theta_{12}^2 + \left( \theta_{22} - 1 \right)^2 \le 1 \;, \qquad
0 \le \theta_{22} \le 2 \;.
\end{equation}
From (\ref{eq:infthetas3}) we learn that
the marginal variance, $\theta^{(\infty)}_{11}$,
depends only on $\tau$, $\delta$ and the correlation of $\Theta$.
Furthermore, if $\Delta$ is not too small,
then if $\delta$ is small (as in Fig.~\ref{fig:stochBifDiag_c}),
$\theta^{(\infty)}_{22} \approx \theta_{22}$,
and if $\tau$ is small (as in Fig.~\ref{fig:stochBifDiag_b}),
$\theta^{(\infty)}_{11}$ may assume only a small range of values.

\section{Invariant densities}
\label{sec:INV}
\setcounter{equation}{0}

\subsection{Near Gaussian densities}
\label{sub:D1}

As described at the beginning of \S\ref{sec:LINEAR},
if (\ref{eq:P}) with $\ee = 0$ has a periodic solution,
$\{ (x_i,y_i) \}_{i = 0,\ldots,n-1}$,
satisfying $x_i \ne 0$ for each $i$,
and no other attractor,
then for small $\ee > 0$
we observe an invariant density that is well-approximated by the sum of $n$ Gaussians
centred about each $(x_i,y_i)$, and scaled by $\frac{1}{n}$.
The size of each Gaussian is $O(\ee)$.
The precise shape of each Gaussian approximation
is given by the results of \S\ref{sub:GAUSS}.
Specifically its mean is $(x_i,y_i)$,
and its covariance matrix, $\Lambda^{(i)}$, is given iteratively by (\ref{eq:Lambdai}),
where $\Lambda^{(0)}$ is given by (\ref{eq:lambdas}).

\begin{figure}[t!]
\begin{center}
\setlength{\unitlength}{1cm}
\begin{picture}(16,4.1)
\put(.4,.3){\includegraphics[height=3.8cm]{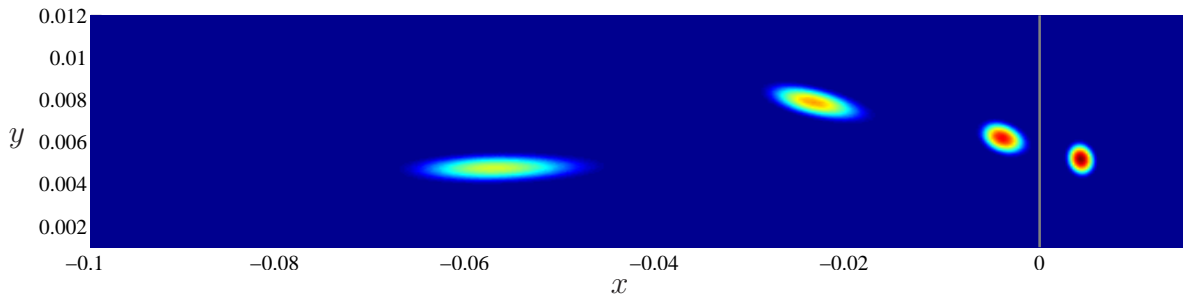}}
\put(8,0){$x$}
\put(0,2){$y$}
\end{picture}
\caption{
The invariant density of the stochastic Nordmark map, (\ref{eq:P}),
using the same parameter values as Fig.~\ref{fig:stochBifDiag_c}
with also $\mu = 0.005$.
Here the noise amplitude is sufficiently small ($\ee = 0.00025$)
that the density is approximately the sum of four Gaussians, scaled by $\frac{1}{4}$.
The value of the density is indicated by color
(dark red -- the maximum value of the density;
dark blue -- zero).
The vertical gray line denotes the switching manifold, $x=0$.
\label{fig:invDensity_a}
}
\end{center}
\end{figure}

An example is shown in Fig.~\ref{fig:invDensity_a}.
This figure, and others below,
shows the invariant probability density function
where the value of the function is indicated by colour
(dark red -- the maximum value of the density;
dark blue -- zero)\removableFootnote{
In order to make the figures clear
I had to make the colour scale nonlinear
in that a wider range of colours is used for smaller values of the density.
For instance, green, which is about half way between dark red
and dark blue, corresponds to a value of the density equal to about $\frac{1}{10}$
of its maximum value (instead of about $\frac{1}{2}$).
This nonlinearity is necessary because, for instance,
roughly speaking, each of the four humps in Fig.~\ref{fig:invDensity_a}
corresponds to a quarter of the iterates,
thus has a volume of $\frac{1}{4}$,
so since the covariance of the right-most hump is smaller than that of the left-most hump,
the peak of the right-most hump is significantly higher than that of the left-most hump,
and so with a linear colour scale the left-most hump
would appear as merely an innocuous blue blur.
}.
Following \cite{Gr05}, we assumed the system is ergodic and approximated invariant densities
by a two-dimensional histogram computed from $10^8$ consecutive iterates of a single orbit.

With the parameter values of Fig.~\ref{fig:invDensity_a},
iterates of (\ref{eq:P}) follow close to the underlying deterministic period-$4$ solution.
Iterates step from the neighbourhood of each point in the order left to right.
While $x<0$, iterates approach the deterministic fixed point of the left half-map, (\ref{eq:xyStarL}), which for $\mu > 0$ lies in the right half-plane and is virtual.
Near the right-most point of the periodic solution, for which $x>0$,
the square-root term takes effect and the next iterate
lies near the left-most point.
The left-most Gaussian has the greatest spread because
it corresponds to iterates that have just experienced the
strongly expanding square-root term.
Moving left to right,
the size of each Gaussian decreases because the left half-map of (\ref{eq:P})
is contracting, and indeed we have the explicit formula (\ref{eq:Lambdai}).
Note that since the left-most Gaussian is most widely spread,
the value of the density at its peak is substantially less than
at the other three peaks.

\subsection{Weakly non-Gaussian densities}
\label{sub:D2}

\begin{figure}[t!]
\begin{center}
\setlength{\unitlength}{1cm}
\begin{picture}(8,4.3)
\put(0,0){\includegraphics[height=4cm]{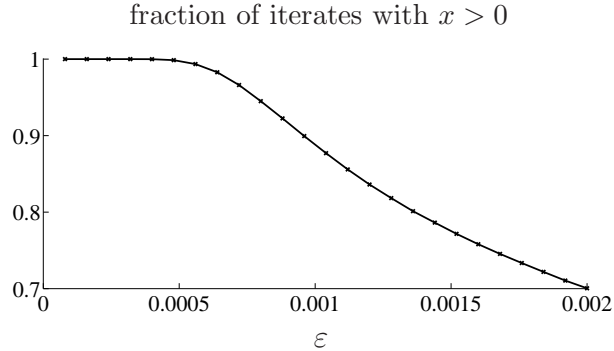}}
\put(4,0){$\ee$}
\put(1.55,4.3){\small fraction of iterates with $x>0$}
\end{picture}
\caption{
The fraction of points of an orbit
(with transients decayed) of (\ref{eq:P}) with $x>0$
that return to the right half-plane in exactly four iterations of (\ref{eq:P}).
The parameter values are the same as in Fig.~\ref{fig:stochBifDiag_c}.
Each data point was computed numerically from $10^7$ iterations of one orbit.
The invariant density when $\ee = 0.00025$ is shown in Fig.~\ref{fig:invDensity_a},
and when $\ee = 0.00075$ in Fig.~\ref{fig:invDensity_b}.
\label{fig:periodFrac}
}
\end{center}
\end{figure}

As the value of $\ee$ is increased,
the above approximation to the invariant density worsens because
with increasing probability iterates fall on the
``wrong'' side of the switching manifold.
For an arbitrary orbit,
let $\mathcal{I}_j$ denote the number of iterations
that the $j^{\rm th}$-point of the orbit satisfying $x>0$
takes to return to the right half-plane.
Fig.~\ref{fig:periodFrac} is a plot of the fraction of $\mathcal{I}_j$
that are equal to $n$,
where $n$ is the period of the deterministic periodic solution, against $\ee$,
using the same parameter values as Fig.~\ref{fig:stochBifDiag_c}, and $\mu = 0.005$.

A key feature of Fig.~\ref{fig:periodFrac} is that,
roughly speaking, this fraction begins to descend rapidly
over a narrow interval of values of $\ee$ (around $\ee = 0.0006$).
Two main factors influence the location of this critical interval.
If the deterministic periodic solution is weakly attracting,
then the matrix, $K$, (\ref{eq:K}), has an eigenvalue near the unit circle,
and thus by (\ref{eq:detImM}) the matrix, $I-M$, is almost singular.
Therefore by (\ref{eq:lambdas}) the invariant density is broad
and so the values of $\ee$ in the critical interval are relatively small.
Second, the values of $\ee$ in the interval are roughly proportional
to the distance of the deterministic periodic solution from the switching manifold.
These observations explain the breakdown of the Gaussian approximation visible in
Figs.~\ref{fig:stochBifDiag_c} and \ref{fig:stochBifDiag_b}.
Near a saddle-node or period-doubling bifurcation
the strength of attraction of the periodic solution is weak and the Gaussian fit is poor.
Near a border-collision bifurcation
the periodic solution is in close proximity to the switching manifold
and again iterations often stray from the Gaussian approximation.

\begin{figure}[t!]
\begin{center}
\setlength{\unitlength}{1cm}
\begin{picture}(16,4.1)
\put(.4,.3){\includegraphics[height=3.8cm]{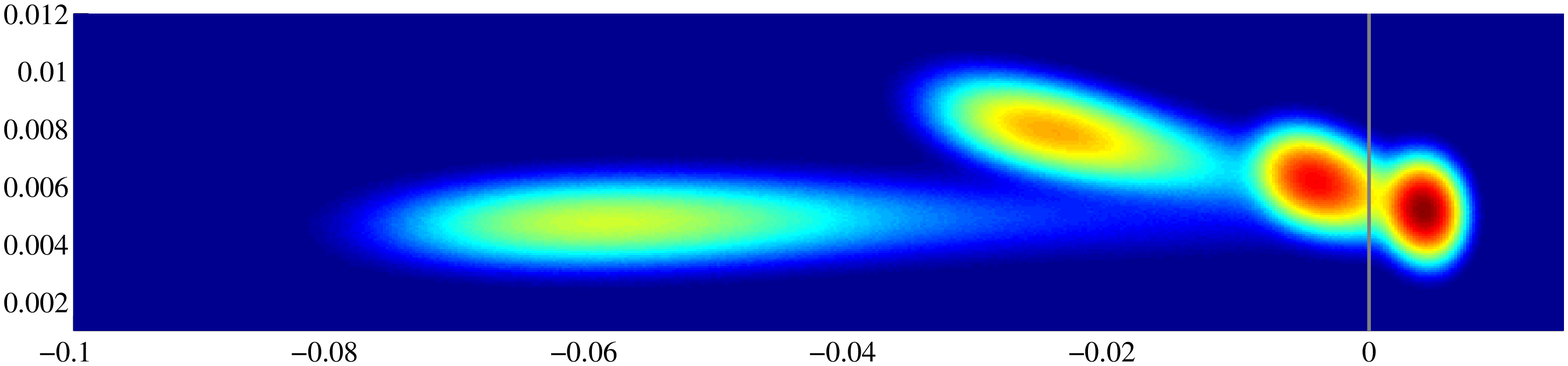}}
\put(8,0){$x$}
\put(0,2){$y$}
\end{picture}
\caption{
The invariant density of (\ref{eq:P})
when $\mu = 0.005$, $\ee = 0.00075$
and the remaining parameter values are as in Fig.~\ref{fig:stochBifDiag_c}.
\label{fig:invDensity_b}
}
\end{center}
\end{figure}

When $\ee = 0.00025$, as in Fig.~\ref{fig:invDensity_a},
iterations follow the deterministic period-$4$ pattern almost exclusively
and the Gaussian approximation is a good fit to the
numerically computed invariant density, as described above.
When $\ee = 0.00075$, as in Fig.~\ref{fig:periodFrac},
iterates follow the period-$4$ pattern about $96\%$ of the time.
In this case the invariant density, shown in Fig.~\ref{fig:invDensity_b},
exhibits four clear peaks near the deterministic solution,
but level curves are relatively non-elliptical indicating that the density
is not as well approximated by Gaussians.
A relatively high fraction of iterates
lie very near the switching manifold.
Due to the square-root singularity, (\ref{eq:P}) maps
a small region just to the right of the switching manifold
to a large region that including points relatively
distant from the deterministic periodic solution,
say around $(x,y) \approx (-0.02,0.005)$.
Consequently the invariant density displays a nonlinear character.

\subsection{Strongly non-Gaussian densities}
\label{sub:D3}

\begin{figure}[t!]
\begin{center}
\setlength{\unitlength}{1cm}
\begin{picture}(16,4.1)
\put(.4,.3){\includegraphics[height=3.8cm]{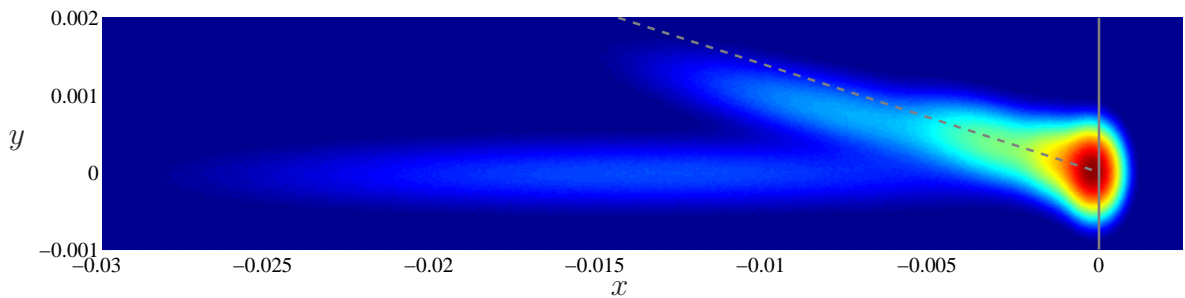}}
\put(8,0){$x$}
\put(0,2){$y$}
\end{picture}
\caption{
The invariant density of (\ref{eq:P}) when $\mu = 0$
and the remaining parameter values are as in Fig.~\ref{fig:stochBifDiag_c}.
The dashed line is the slow manifold for the fixed point of the 
left half-map of (\ref{eq:P}).
Since the left half-map is linear, this manifold
coincides with the eigenspace for the stability multiplier of greatest magnitude.
\label{fig:invDensity_c}
}
\end{center}
\end{figure}

\begin{figure}[t!]
\begin{center}
\setlength{\unitlength}{1cm}
\begin{picture}(10.2,5)
\put(.3,0){\includegraphics[height=5cm]{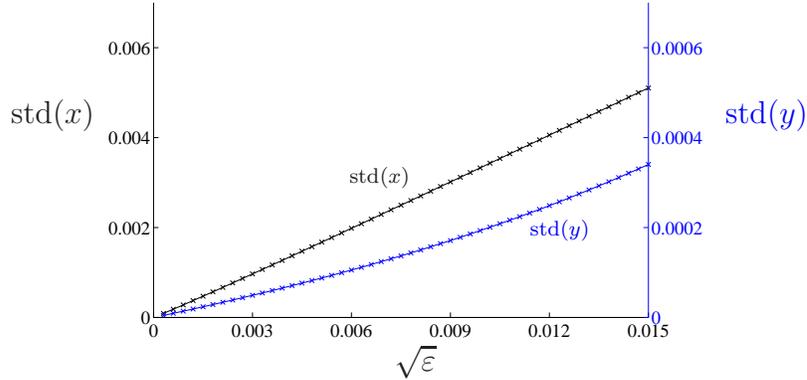}}
\put(5.1,0){$\sqrt{\ee}$}
\put(0,3.3){${\rm std}(x)$}
\put(9.5,3.3){\color{blue} ${\rm std}(y)$}
\put(4.5,2.5){\scriptsize ${\rm std}(x)$}
\put(6.9,1.8){\scriptsize \color{blue} ${\rm std}(y)$}
\end{picture}
\caption{
Standard deviations of the marginal densities of the invariant density
(shown in Fig.~\ref{fig:invDensity_c} for $\ee = 0.00025$)
of (\ref{eq:P}) for different values of $\ee$ when $\mu = 0$
and the remaining parameter values are as in Fig.~\ref{fig:stochBifDiag_c}.
Each data point was computed numerically from $10^7$ iterations of one orbit
with transients decayed.
\label{fig:deviationEe}
}
\end{center}
\end{figure}

If $\ee$ is large relative to
the distance of the deterministic attractor from the switching manifold,
or if this distance is zero,
the invariant density takes a strongly non-Gaussian form.
As an example, Fig.~\ref{fig:invDensity_c}
shows an invariant density of (\ref{eq:P}) at the grazing bifurcation ($\mu = 0$).
Since we have a thorough understanding of the behaviour of map when $\ee = 0$,
the shape of the density can be explained heuristically.
Due to the square-root singularity,
the map (\ref{eq:P}) flings points with $x>0$
(and $y < \chi \sqrt{x} - \tau x$) relatively far into the left half-plane.
While $x<0$, subsequent iterates are governed by the left half-map
and slowly step towards the fixed point of this map, which, when $\mu = 0$,
is located at the origin.
For Fig.~\ref{fig:invDensity_c}, the stability multipliers associated with
the fixed point are real and distinct (approximately $0.14$ and $0.36$).
Thus iterates first step quickly to the eigenspace corresponding to the
multiplier of greater magnitude ($0.36$)
(indicated in Fig.~\ref{fig:invDensity_c}),
then head into the origin.
This explains the angle at which the invariant density protrudes from $x=0$.
Since noise is present,
iterates eventually pass over the switching manifold,
and the process repeats.
Since the noise amplitude is $O(\ee)$, iterates with $x>0$ are $O(\ee)$.
Due to the square-root singularity,
the size of the invariant density is therefore $O(\sqrt{\ee})$, Fig.~\ref{fig:deviationEe}.
In addition, when $\mu = 0$ the average number of iterations taken for points of an orbit to
return to the right half-space decreases as $\ee$ is increased.

\subsection{Densities resulting from coexisting attractors}
\label{sub:D4}

\begin{figure}[t!]
\begin{center}
\setlength{\unitlength}{1cm}
\begin{picture}(10,5.1)
\put(.3,.3){\includegraphics[height=4.8cm]{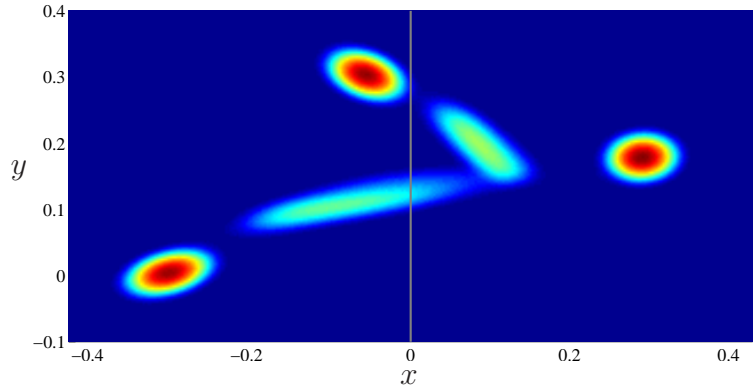}}
\put(5.18,0){$x$}
\put(0,2.8){$y$}
\end{picture}
\caption{
The invariant density of (\ref{eq:P}) when $\mu = 0.15$
and the remaining parameter values are as in Fig.~\ref{fig:stochBifDiag_b}.
\label{fig:invDensity_d}
}
\end{center}
\end{figure}

If, in the absence of noise (\ref{eq:P}) has coexisting attractors,
with noise we expect there to be a single invariant density involving peaks
near both deterministic attractors.
Such a density is shown in Fig.~\ref{fig:invDensity_d}.
The relative weighting of the density about each deterministic attractor
depends on the relative strengths of attraction
and the relative sizes of the basins of attraction.
Orbits may remain in the neighbourhood of one attractor
for a long period of time before sufficiently extreme noise
pushes the orbit into the neighbourhood of a different attractor.
Similar behaviour has been described for smooth maps \cite{KnWe89}.

\section{Conclusions}
\label{sec:CONC}
\setcounter{equation}{0}

In this paper we studied
the Nordmark map with additive, Gaussian noise, (\ref{eq:P}),
by analyzing two-dimensional invariant densities that
were computed numerically and in certain situations approximated
by a sum of scaled Gaussian probability density functions.

In \S\ref{sec:DERIV} we derived (\ref{eq:P})
from a general three-dimensional stochastic differential equation, (\ref{eq:SDE}),
that, in the absence of noise, exhibits a regular grazing bifurcation.
The map (\ref{eq:P}) is the composition of a smooth map encapsulating global dynamics,
(the global map, (\ref{eq:G}))
and a piecewise-smooth local map (the discontinuity map, (\ref{eq:D}))
that represents the correction to the global map induced by the dynamics for $u > 0$.
The noise term in (\ref{eq:G}) was calculated
using standard stochastic asymptotics theory, \S\ref{sub:ADD}.
Near grazing the discontinuity map corresponds to dynamics occurring over a much shorter time period
than the global map and for this reason
the effect of randomness in the discontinuity map may be neglected.
Moreover, with this assumption we derived the covariance matrix, $\Theta$, for the noise
in (\ref{eq:P}) explicitly in terms of the components of the equations of motion, (\ref{eq:Theta}).

The map (\ref{eq:P}) cannot be transformed in a manner that puts
$\Theta$ into a canonical form without disrupting the deterministic normal form.
In \S\ref{sub:APPROX} and \S\ref{sub:THETA} we showed that with certain reasonable assumptions,
the shape of $\Theta$ has little effect on the resulting invariant density.

Numerical simulations of (\ref{eq:P}) revealed
invariant densities that roughly conformed to one of four descriptions.
If, in the absence of noise ($\ee = 0$) there exists an attracting period-$n$ solution 
with no points on the switching manifold,
then for sufficiently small $\ee > 0$
the invariant density is approximately the sum of $n$ Gaussians
scaled by $\frac{1}{n}$, \S\ref{sub:D1}.
The associated covariance matrices were computed in \S\ref{sub:GAUSS}.
In the same situation with larger values of $\ee$,
iterates fall near the switching manifold with higher probability.
Since the square-root singularity in (\ref{eq:P}) induces
a strongly nonlinear stretching of phase-space near the switching manifold,
the Gaussian approximation worsens as $\ee$ is increased, \S\ref{sub:D2}.
With still larger $\ee$,
the invariant density may be highly irregular
but of a shape amenable to a simple explanation, \S\ref{sub:D3}.
When attractors coexist in the absence of noise,
for $\ee > 0$ the invariant density is comprised of pieces about each deterministic
attractor, \S\ref{sub:D4}.
The rough magnitude of $\ee$ for which the effect of the square-root singularity becomes significant, 
is influenced by the strength of attraction of the deterministic period-$n$ solution
and the proximity of the solution to the switching manifold.

Given the above descriptions of the numerically computed invariant density
at various fixed values of $\mu$,
we can now explain the dominate features of
Figs.~\ref{fig:stochBifDiag_c} and \ref{fig:stochBifDiag_b}.
We note that transitions in Figs.~\ref{fig:stochBifDiag_c} and \ref{fig:stochBifDiag_b}
are {\em P-bifurcations} \cite{Ar98}
at which the invariant density changes character fundamentally.
The sequences of transitions have not been described before
as the underlying deterministic system is piecewise-smooth.

In Fig.~\ref{fig:stochBifDiag_c},
for $\mu < -0.001$, say,
the invariant density is an approximate Gaussian centred at
the fixed point, (\ref{eq:xyStarL}).
Note that the width of the invariant density about the non-impacting solution
is practically constant leading up to the grazing bifurcation,
which is not the typical case for smooth maps \cite{CrFa82,LiLu86}.
For $\mu \le 0$, with $|\mu| \le O(\ee)$,
the fixed point is close enough to the switching manifold that
iterates repeatedly lie in the right half-plane.
Due to the square-root singularity,
the invariant density is non-Gaussian and $O(\sqrt{\ee})$ in size.

For small $\mu > 0$, there is an irregular invariant density of size $O(\sqrt{\ee})$;
roughly periodic behaviour at high periods is not seen.
Over a wide interval of intermediate values of $\mu$
the invariant density is well approximated by four Gaussians.
As $\mu$ increases across this interval, the approximation
improves because the strength of attraction of
the deterministic period-$4$ solution increases with $\mu$.
The four-piece invariant density transforms to
a widely spread density as the deterministic period-$4$ solution
collides with the switching manifold in a border-collision bifurcation
at $\mu \approx 0.014$.
In the absence of noise, just beyond this bifurcation, there is a chaotic attractor.
With noise the invariant density has a complicated wide-spread structure.
For larger $\mu$ the invariant density is well-approximated by three Gaussians.

Fig.~\ref{fig:stochBifDiag_b} displays much the same behaviour,
but in addition iterates traverse
neighbourhoods of coexisting deterministic attractors.
Note that in some cases iterates only appear to
surround a single attractor (e.g.~for $\mu = 0.2$).
This occurs when the attractor is strongly attracting and its basin
of attraction is large, relative to the noise amplitude,
and the number of iterates of the numerical simulation is relatively small.

It remains to investigate grazing bifurcations of vibro-impacting systems
by studying a model with noise applied only to impact events.
The formulation (\ref{eq:P}) does not allow for an investigation into
rapidly occurring, noise-induced multiple impacts.
Many other bifurcation scenarios in piecewise-smooth systems
await a detailed investigation in the presence of noise, such as sliding bifurcations.





\appendix
\section{The stochastic Nordmark map for the linear oscillator}
\label{sec:APP}
\setcounter{equation}{0}

Here we derive the parameters,
$\tau$, $\delta$, $\chi$, and $\Theta$,
of the stochastic Nordmark map, (\ref{eq:P}),
for the impact oscillator described in \S\ref{sub:IMPACT}.
We first let
\begin{equation}
\bu = \left[ \begin{array}{c}
u \\ v \\ w
\end{array} \right] = \left[ \begin{array}{c}
u \\ \dot{u} \\ (t-t_{\rm graz}) {\rm ~mod~} 2 \pi
\end{array} \right] \;, \qquad
\eta = F - F_{\rm graz} \;,
\end{equation}
where
\begin{equation}
t_{\rm graz} = \tan^{-1} \left( \frac{b_{\rm osc}}{k_{\rm osc}-1} \right) \;,
\label{eq:tGraz}
\end{equation}
is the time in $\left( \frac{\pi}{2}, \pi \right)$
at which the impact occurs for the grazing periodic orbit,
so that the system takes the general form (\ref{eq:SDE})
and satisfies the related assumptions
(except phase space is $\mathbb{R}^2 \times \mathbb{S}^1$ instead of
$\mathbb{R}^3$, but this is inconsequential).
From the deterministic solution to (\ref{eq:impactOsc})\removableFootnote{
The deterministic solution to (\ref{eq:impactOsc}) is
\begin{equation}
\left[ \begin{array}{c} u_{\rm det}(t) \\ \dot{u}_{\rm det}(t) \end{array} \right] =
\left[ \begin{array}{c}
-1 + \kappa_1 \cos(t) + \kappa_2 \sin(t) \\
-\kappa_1 \sin(t) + \kappa_2 \cos(t)
\end{array} \right] +
{\rm e}^{J(t-t_0)}
\left[ \begin{array}{c}
u_0 + 1 - \kappa_1 \cos(t_0) - \kappa_2 \sin(t_0) \\
\dot{u}_0 + \kappa_1 \sin(t_0) - \kappa_2 \cos(t_0)
\end{array} \right] \;,
\label{eq:detSoln}
\end{equation}
where
\begin{equation}
\left[ \begin{array}{c} \kappa_1 \\ \kappa_2 \end{array} \right] =
\frac{F}{b_{\rm osc}^2 + (1-k_{\rm osc})^2}
\left[ \begin{array}{c} k_{\rm osc}-1 \\ b_{\rm osc} \end{array} \right] \;.
\end{equation}
},
we find that the global map, $G$, is given by (\ref{eq:G}) with
\begin{equation}
\hat{A} = {\rm e}^{2 \pi J} \;, \qquad
\hat{b} = \frac{1}{F_{\rm graz}} \left[ \begin{array}{c}
1 - \hat{a}_{11} \\ -\hat{a}_{21}
\end{array} \right] \;,
\end{equation}
where
\begin{equation}
J = \left[ \begin{array}{cc} 0 & 1 \\ -k_{\rm osc} & -b_{\rm osc} \end{array} \right] \;.
\label{eq:J}
\end{equation}
Also
\begin{equation}
c = \frac{2 \sqrt{2} k_{\rm supp} d}{1 + k_{\rm supp} d} \;,
\end{equation}
in the discontinuity map, $D$, (\ref{eq:D}).
From the coordinate change (\ref{eq:coordinateChange}), we obtain
\begin{equation}
\tau = 2 {\rm e}^{2 \pi \alpha} \cos(2 \pi \beta) \;, \qquad
\delta = {\rm e}^{4 \pi \alpha} \;, \qquad
\chi = {\rm sgn}(\hat{a}_{12} c) \;,
\end{equation}
where $\alpha \pm {\rm i} \beta$ are the eigenvalues of $J$, and specifically,
\begin{equation}
\alpha = -\frac{b_{\rm osc}}{2} \;, \qquad
\beta = \sqrt{k_{\rm osc} - \frac{b_{\rm osc}^2}{4}} \;.
\end{equation}

It remains to use (\ref{eq:Theta}) to compute the covariance matrix,
$\Theta$, of the noise, $\xi$.
Here we evaluate all terms of (\ref{eq:Theta}), for (\ref{eq:impactOsc}).
In the form (\ref{eq:SDE}), we have
\begin{eqnarray}
f^{(L)}(\bu;\eta) &=& \left[ \begin{array}{c}
v \\ -k_{\rm osc} (u+1) - b_{\rm osc} v + (\eta + F_{\rm graz}) \cos(w + t_{\rm graz}) \\ 1
\end{array} \right] \;, \\
B(\bu;\eta) &=& \left[ \begin{array}{ccc}
0 & 0 & 0 \\
0 & 1 & 0 \\
0 & 0 & 0
\end{array} \right] \;,
\end{eqnarray}
so by (\ref{eq:FGraz}) and (\ref{eq:tGraz}),
$f^{(L)} \left( [0,0,0]^{\sf T} ; 0 \right) = [0,-\gamma_L,\zeta_L]^{\sf T}
= [0,-1,1]^{\sf T}$.
Then (\ref{eq:H}) yields
\begin{equation}
H(s,t) = \left[ \begin{array}{ccc}
0 & {\rm e}^{J(t-s)}_{12} & 0 \\
0 & {\rm e}^{J(t-s)}_{22} & 0 \\
0 & 0 & 0
\end{array} \right] \;,
\label{eq:HimpactOsc}
\end{equation}
where ${\rm e}^{J(t-s)}_{ij}$ denotes the $(i,j)$-element of ${\rm e}^{J(t-s)}$.
Finally, the evaluation of (\ref{eq:Omega}) at $T^{(0)} = 2 \pi$,
using (\ref{eq:J}) and (\ref{eq:HimpactOsc}), gives
\begin{eqnarray}
\omega_{11} &=& \frac{1}{4 \alpha \beta^2} ({\rm e}^{2 \alpha t} - 1)
- \frac{\alpha}{4 k_{\rm osc} \beta^2} ({\rm e}^{2 \alpha t} \cos(2 \beta t) - 1)
- \frac{1}{4 k_{\rm osc} \beta} {\rm e}^{2 \alpha t} \sin(2 \beta t) \;, \\
\omega_{12} &=& -\frac{1}{4 \beta^2} {\rm e}^{2 \alpha t} (\cos(2 \beta t) - 1) \;, \\
\omega_{22} &=& \frac{k_{\rm osc}}{4 \alpha \beta^2} ({\rm e}^{2 \alpha t} - 1)
- \frac{\alpha}{4 \beta^2} ({\rm e}^{2 \alpha t} \cos(2 \beta t) - 1)
+ \frac{1}{4 \beta} {\rm e}^{2 \alpha t} \sin(2 \beta t) \;,
\end{eqnarray}
and $\omega_{13} = \omega_{23} = \omega_{33} = 0$\removableFootnote{
These three values are zero because $w$ is deterministic.
}.


\begin{thebibliography}{10}

\bibitem{WiDe00}
M.~Wiercigroch and B.~De~Kraker, editors.
\newblock {\em Applied Nonlinear Dynamics and Chaos of Mechanical Systems with
  Discontinuities.}, Singapore, 2000. World Scientific.

\bibitem{Br99}
B.~Brogliato.
\newblock {\em Nonsmooth Mechanics: Models, Dynamics and Control.}
\newblock Springer-Verlag, New York, 1999.

\bibitem{BlCz99}
B.~Blazejczyk-Okolewska, K.~Czolczynski, T.~Kapitaniak, and J.~Wojewoda.
\newblock {\em Chaotic Mechanics in Systems with Impacts and Friction}.
\newblock World Scientific, Singapore, 1999.

\bibitem{Ib09}
R.A. Ibrahim.
\newblock {\em Vibro-Impact Dynamics.}, volume~43 of {\em Lecture Notes in
  Applied and Computational Mechanics.}
\newblock Springer, New York, 2009.

\bibitem{DaZh07}
H.~Dankowicz, X.~Zhao, and S.~Misra.
\newblock Near-grazing dynamics in tapping-mode atomic-force microscopy.
\newblock {\em Int. J. Non-Linear Mech.}, 42(4):697--709, 2007.

\bibitem{RaMe08}
A.~Raman, J.~Melcher, and R.~Tung.
\newblock Cantilever dynamics in atomic force microscopy.
\newblock {\em Nanotoday}, 3(1-2):20--27, 2008.

\bibitem{ThNa00}
S.~Theodossiades and S.~Natsiavas.
\newblock Non-linear dynamics of gear-pair systems with periodic stiffness and
  backlash.
\newblock {\em J. Sound Vib.}, 229(2):287--310, 2000.

\bibitem{HaWi07}
C.K. Halse, R.E. Wilson, M.~di~Bernardo, and M.E. Homer.
\newblock Coexisting solutions and bifurcations in mechanical oscillations with
  backlash.
\newblock {\em J. Sound Vib.}, 305:854--885, 2007.

\bibitem{Gr88}
I.~Grabec.
\newblock Chaotic dynamics of the cutting process.
\newblock {\em Int. J. Mach. Tools Manufact.}, 28(1):19--32, 1988.

\bibitem{Wi97}
M.~Wiercigroch.
\newblock Chaotic vibration of a simple model of the machine tool-cutting
  process system.
\newblock {\em J. Vib. Acoust.}, 119(3):468--475, 1997.

\bibitem{PaLi92}
M.P. Pa\"{\i}doussis and G.X. Li.
\newblock Cross-flow-induced chaotic vibrations of heat-exchanger tubes
  impacting on loose supports.
\newblock {\em J. Sound Vib.}, 152(2):305--326, 1992.

\bibitem{DeFr99}
J.M. de~Bedout, M.A. Franchek, and A.K. Bajaj.
\newblock Robust control of chaotic vibrations for impacting heat exchanger
  tubes in crossflow.
\newblock {\em J. Sound Vib.}, 227(1):183--204, 1999.

\bibitem{DiBu08}
M.~di~Bernardo, C.J. Budd, A.R. Champneys, and P.~Kowalczyk.
\newblock {\em Piecewise-smooth Dynamical Systems. Theory and Applications.}
\newblock Springer-Verlag, New York, 2008.

\bibitem{LiGa04}
B.~Lindner, J.~Garcia-Ojalvo, A.~Neiman, and L.~Schimansky-Geier.
\newblock Effects of noise in excitable systems.
\newblock {\em Phys. Reports}, 392:321--424, 2004.

\bibitem{PiKu97}
A.S. Pikovsky and J.~Kurths.
\newblock Coherence resonance in a noise-driven excitable system.
\newblock {\em Phys. Rev. Lett.}, 78(5):775--778, 1997.

\bibitem{MuVa05}
C.B. Muratov, E.~Vanden-Eijnden, and E.~Weinan.
\newblock Self-induced stochastic resonance in excitable systems.
\newblock {\em Phys. D}, 210:227--240, 2005.

\bibitem{DeGu12}
M.~Desroches, J.~Guckenheimer, B.~Krauskopf, C.~Kuehn, H.M. Osinga, and
  M.~Wechselberger.
\newblock Mixed-mode oscillations with multiple time scales.
\newblock {\em SIAM Rev.}, 54(2):211--288, 2012.

\bibitem{Ra95}
S.~Rajasekar.
\newblock Controlling of chaotic motion by chaos and noise signals in a
  logistic map and a {B}onhoeffer--van der {P}ol oscillator.
\newblock {\em Phys. Rev. E}, 51(1):775--778, 1995.

\bibitem{HoLe06}
W.~Horsthemke and R.~Lefever.
\newblock {\em Noise-Induced Transitions: {T}heory and Applications in Physics,
  Chemistry, and Biology.}
\newblock Springer, New York, 2006.

\bibitem{Da98}
A.~Daffertshofer.
\newblock Effects of noise on the phase dynamics of nonlinear oscillators.
\newblock {\em Phys. Rev. E}, 58(1):327--338, 1998.

\bibitem{ArBl99}
L.~Arnold, G.~Bleckert, and K.R. Schenk-Hopp\'{e}.
\newblock The stochastic {B}russelator: {P}arametric noise destroys {H}opf
  bifurcation.
\newblock In H.~Crauel and M.~Gundlach, editors, {\em Stochastic Dynamics.},
  pages 71--92, 1999.

\bibitem{ZhLu93}
W.Q. Zhu, M.Q. Lu, and Q.T. Wu.
\newblock Stochastic jump and bifurcation of a {D}uffing oscillator under
  narrow-band excitation.
\newblock {\em J. Sound Vib.}, 165(2):285--304, 1993.

\bibitem{Sc96}
K.R. Schenk-Hopp\'{e}.
\newblock Bifurcation scenarios of the noisy {D}uffing-van der {P}ol
  oscillator.
\newblock {\em Nonlinear Dynamics}, 11:255--274, 1996.

\bibitem{GrHo12}
T.C.L. Griffin and S.J. Hogan.
\newblock The effects of noise on a piecewise linear map.
\newblock In preparation., 2012.

\bibitem{Gr05}
T.C.L. Griffin.
\newblock {\em Dynamics of Stochastic Nonsmooth Systems.}
\newblock PhD thesis, University of Bristol, 2005.

\bibitem{GrHo05}
T.~Griffin and S.~Hogan.
\newblock Dynamics of discontinuous systems with imperfections and noise.
\newblock In G.~Rega and F.~Vestroni, editors, {\em IUTAM Symposium on Chaotic
  Dynamics and Control of Systems and Processes in Mechanics.}, pages 275--285.
  Springer, 2005.

\bibitem{GrHo12b}
T.C.L. Griffin, S.J. Hogan, G.~Olivar, and V.~Moreno.
\newblock The effect of noise on the dynamics of {DC/DC} converters.
\newblock In preparation., 2012.

\bibitem{Wa98}
R.~Wackerbauer.
\newblock Noise-induced stabilization of one-dimensional discontinuous maps.
\newblock {\em Phys. Rev. E}, 58(3):3036--3044, 1998.

\bibitem{DiMe79}
M.F. Dimentberg and A.I. Menyailov.
\newblock Response of a single-mass vibroimpact system to white-noise random
  excitation.
\newblock {\em Z. Angew. Math. Mech.}, 59(12):709--716, 1979.

\bibitem{Zh76}
V.F. Zhuravlev.
\newblock A method for analyzing vibration-impact systems by means of special
  functions.
\newblock {\em Mech. Solids}, 11:23--27, 1976.

\bibitem{RoSp86}
J.B. Roberts and P.D. Spanos.
\newblock Stochastic averaging: {A}n approximate method of solving random
  vibration problems.
\newblock {\em Int. J. Non-Linear Mechanics}, 21(2):111--134, 1986.

\bibitem{Kh66}
R.Z. Khas'minskii.
\newblock A limit theorem for the solutions of differential equations with
  random right-hand sides.
\newblock {\em Theory Probab. Appl.}, 11:390--405, 1966.

\bibitem{AnAs02}
V.S. Anishchenko, V.V. Astakhov, A.B. Neiman, T.E. Vadivasova, and
  L.~Schimansky-Geier.
\newblock {\em Nonlinear Dynamics of Chaotic and Stochastic Systems. Tutorial
  and Modern Developments.}
\newblock Springer, New York, 2002.

\bibitem{FoBr96}
M.~Fogli, P.~Bressolette, and P.~Bernard.
\newblock The dynamics of a stochastic oscillator with impacts.
\newblock {\em Eur. J. Mech. A-Solids}, 15(2):213--241, 1996.

\bibitem{DiIo04}
M.F. Dimentberg and D.V. Iourtchenko.
\newblock Random vibrations with impacts: {A} review.
\newblock {\em Nonlinear Dyn.}, 36:229--254, 2004.

\bibitem{SrPa05}
N~Sri~Namachchivaya and J.H. Park.
\newblock Stochastic dynamics of impact oscillators.
\newblock {\em J. Appl. Mech. Trans. ASME}, 72(6):862--870, 2005.

\bibitem{DiIo05}
M.F. Dimentberg and D.V. Iourtchenko.
\newblock Stochastic and/or chaotic response of a vibration system to
  imperfectly periodic sinusoidal excitation.
\newblock {\em Int. J. Bifurcation Chaos}, 15(6):2057--2061, 2005.

\bibitem{FePf98}
Q.~Feng and F.~Pfeiffer.
\newblock Stochastic model on a rattling system.
\newblock {\em J. Sound Vib.}, 215(3):439--453, 1998.

\bibitem{PfKu90}
F.~Pfeiffer and A.~Kunert.
\newblock Rattling models from deterministic to stochastic processes.
\newblock {\em Nonlinear Dyn.}, 1:63--74, 1990.

\bibitem{DiBu01}
M.~di~Bernardo, C.J. Budd, and A.R. Champneys.
\newblock Normal form maps for grazing bifurcations in $n$-dimensional
  piecewise-smooth dynamical systems.
\newblock {\em Phys. D}, 160:222--254, 2001.

\bibitem{No91}
A.B. Nordmark.
\newblock Non-periodic motion caused by grazing incidence in impact
  oscillators.
\newblock {\em J. Sound Vib.}, 2:279--297, 1991.

\bibitem{No97}
A.B. Nordmark.
\newblock Universal limit mapping in grazing bifurcations.
\newblock {\em Phys. Rev. E}, 55(1):266--270, 1997.

\bibitem{No01}
A.B. Nordmark.
\newblock Existence of periodic orbits in grazing bifurcations of impacting
  mechanical oscillators.
\newblock {\em Nonlinearity}, 14:1517--1542, 2001.

\bibitem{CrFa82}
J.P. Crutchfield, J.D. Farmer, and B.A. Huberman.
\newblock Fluctuations and simple chaotic dynamics.
\newblock {\em Phys. Rep.}, 92(2):45--82, 1982.

\bibitem{MaHa81}
G.~Mayer-Kress and H.~Haken.
\newblock The influence of noise on the logistic model.
\newblock {\em J. Stat. Phys.}, 26(1):149--171, 1981.

\bibitem{OeHi97}
M.~Oestreich, N.~Hinrichs, K.~Popp, and C.J. Budd.
\newblock Analytical and experimental investigation of an impact oscillator.
\newblock In {\em Proceedings of the ASME 16th Biennal Conf. on Mech.
  Vibrations and Noise.}, pages 1--11, 1997.

\bibitem{DaNo00}
H.~Dankowicz and A.B. Nordmark.
\newblock On the origin and bifurcations of stick-slip oscillations.
\newblock {\em Phys. D}, 136:280--302, 2000.

\bibitem{DiBu01c}
M.~di~Bernardo, C.J. Budd, and A.R. Champneys.
\newblock Corner collision implies border-collision bifurcation.
\newblock {\em Phys. D}, 154:171--194, 2001.

\bibitem{FrNo97}
M.H. Fredriksson and A.B. Nordmark.
\newblock Bifurcations caused by grazing incidence in many degrees of freedom
  impact oscillators.
\newblock {\em Proc. R. Soc. A}, 453:1261--1276, 1997.

\bibitem{FrNo00}
M.H. Fredriksson and A.B. Nordmark.
\newblock On normal form calculation in impact oscillators.
\newblock {\em Proc. R. Soc. A}, 456:315--329, 2000.

\bibitem{MoDe01}
J.~Molenaar, J.G. de~Weger, and W.~van~de Water.
\newblock Mappings of grazing-impact oscillators.
\newblock {\em Nonlinearity}, 14:301--321, 2001.

\bibitem{PrSh98}
Yu.V. Prokhorov and A.N. Shiryaev, editors.
\newblock {\em Probability Theory III: Stochastic Calculus.}
\newblock Springer, New York, 1998.

\bibitem{Ka90}
T.~Kapitaniak.
\newblock {\em Chaos in Systems with Noise.}
\newblock World Scientific, Singapore, 1990.

\bibitem{WeKn90}
J.B. Weiss and E.~Knobloch.
\newblock A stochastic return map for stochastic differential equations.
\newblock {\em J. Stat. Phys.}, 58(5-6):863--883, 1990.

\bibitem{FrWe84}
M.I. Freidlin and A.D. Wentzell.
\newblock {\em Random Perturbations of Dynamical Systems.}
\newblock Springer-Verlag, New York, 1984.

\bibitem{Sc10}
Z.~Schuss.
\newblock {\em Theory and Applications of Stochastic Processes.}
\newblock Springer, New York, 2010.

\bibitem{GrVa99}
J.~Grasman and O.A. van Herwaarden.
\newblock {\em Asymptotic Methods for the Fokker-Planck Equation and the Exit
  Problem in Applications.}
\newblock Springer, New York, 1999.

\bibitem{MaIn08}
Y.~Ma, J.~Ing, S.~Banerjee, M.~Wiercigroch, and E.~Pavlovskaia.
\newblock The nature of the normal form map for soft impacting systems.
\newblock {\em Int. J. Nonlinear Mech.}, 43:504--513, 2008.

\bibitem{InPa08b}
J.~Ing, E.~Pavlovskaia, M.~Wiercigroch, and S.~Banerjee.
\newblock Experimental study of impact oscillator with one-sided elastic
  constraint.
\newblock {\em Phil. Trans. R. Soc. A}, 366:679--704, 2008.

\bibitem{InPa06}
J.~Ing, E.~Pavlovskaia, and M.~Wiercigroch.
\newblock Dynamics of a nearly symmetrical piecewise linear oscillator close to
  grazing incidence: {M}odelling and experimental verification.
\newblock {\em Nonlinear Dyn.}, 46:225--238, 2006.

\bibitem{ChOt94}
W.~Chin, E.~Ott, H.E. Nusse, and C.~Grebogi.
\newblock Grazing bifurcations in impact oscillators.
\newblock {\em Phys. Rev. E}, 50(6):4427--4450, 1994.

\bibitem{So02}
T.~S\"{o}derstr\"{o}m.
\newblock {\em Discrete-time Stochastic Systems.}
\newblock Springer, New York, 2002.

\bibitem{LaMa94}
A.~Lasota and M.C. Mackey.
\newblock {\em Chaos, Fractals, and Noise. Stochastic Aspects of Dynamics.}
\newblock Springer-Verlag, New York, 1994.

\bibitem{KnWe89}
E.~Knobloch and J.B. Weiss.
\newblock Effect of noise on discrete dynamical systems with multiple
  attractors.
\newblock In M.F. McClintock and P.V.E. Moss, editors, {\em Noise in Nonlinear
  Dynamical Systems. Theory of noise induced processes in special
  applications.}, volume~2. Cambridge University Press, New York, 1989.

\bibitem{Ar98}
L.~Arnold.
\newblock {\em Random Dynamical Systems.}
\newblock Springer, New York, 1998.

\bibitem{LiLu86}
S.J. Linz and M.~L\"{u}cke.
\newblock Effect of additive and multiplicative noise on the first bifurcations
  of the logistic model.
\newblock {\em Phys. Rev. A}, 33(4):2694--2704, 1986.

\end{thebibliography}

\end{document}